\documentclass[12pt]{amsart}
\usepackage{amssymb}
\usepackage{enumerate}
\usepackage{amsmath}
\usepackage{graphicx,color}

\renewcommand{\le}{\leqslant}
\renewcommand{\ge}{\geqslant}

\newcommand{\setm}{\setminus}
\newcommand{\ope}{\operatorname}
\newcommand{\wt}{\widetilde}
\newcommand{\df}{\buildrel\rm{def}\over=}

\newcommand{\eps}{\varepsilon}
\newcommand{\D}{\mathcal D}
\newcommand{\R}{\mathbb R}
\newcommand{\B}{\mathcal B}
\newcommand{\A}{\mathcal A}

\newcommand{\ci}[1]{_{{}_{\scriptstyle #1}}}

\begin{document}

\title[Bellman function: a case study]{On a Bellman function associated with \\
the Chang--Wilson--Wolff theorem:\\ a case study}

\author{Fedor Nazarov}
\address{Kent State University}
\email{nazarov@math.kent.edu}

\author{Vasily Vasyunin}
\address{St.-Petersburg department\\
of V.~A.~Steklov Mathematical institute\\ 
of the Russian Academy of Sciences;\\
St.-Petersburg State University}
\email{vasyunin@pdmi.ras.ru}

\author{Alexander Volberg}
\address{Michigan State University}
\email{volberg@math.msu.edu}

\begin{abstract}
In this paper we estimate the tail of distribution (i.\,e., the measure of the set $\{f\ge x\}$) for those
functions $f$ whose dyadic square function is bounded by a given constant. In particular we get a bit better
estimate than the estimate following from the Chang--Wilson--Wolf theorem. In the paper we investigate the
Bellman function corresponding to the problem. A curious structure of this function is found: it has jumps
of the first derivative at a dense subset of interval $[0,1]$ (where it is calculated exactly), but it is
of $C^\infty$-class for $x>\sqrt3$ (where it is calculated up to a multiplicative constant).

An unusual feature of the paper consists in the usage of computer calculations in the proof. Nevertheless, all the
proofs are quite rigorous, since only the integer arithmetic was assigned to computer.
\end{abstract}

\keywords{Bellman function, square function, Chang--Wilson--Wolf theorem, supersolutions, distribution function}

\thanks{The third author was supported by NSF grants DMS~1600065 and DMS~19000286}

\maketitle

\setcounter{section}{-1}
\section{Level 0: What to keep in mind when reading this paper}

\subsection{Organization of the paper}

Since this paper is quite technical in some places, we decided to write the text not
in the usual ``linear'' manner where each statement is immediately followed by its proof
and each proof contains all the needed auxiliary statements but rather in a ``tree-like'' manner
where the top level is occupied by just the statements of the main results, the second level
is occupied by the statements of the auxiliary results and the proofs of the main results without
some technical details, the third level is occupied by the technical details missing in
the second level and so on until we reach the last fifth level, which contains the proof of
some specific numerical inequality needed before. So, the reader who wants only to get a general 
impression of what has been done in this article can read just Level~1; the reader who wants, in 
addition, to get a general idea of how everything is proved can stop reading at Level~2, and so on. 

Such a structure means that at each level we will freely use the results from the next levels
and the notation from the previous ones. Within each level we do employ the usual linear structure.

\subsection{Warning about computer assisted proofs}

Many of our proofs of various ``elementary inequalities''
are computer assisted. On the other hand, our standards for using computers in the proofs
are quite strict: we allow only algebraic symbolic manipulation of rational functions and basic integer
 arithmetic. All our computations were done using the Mathematica program by Wolfram Research run
on the Windows XP platform. We believe that there were no bugs in the software that could affect our 
results but, of course, the reader is welcome to check the computations using different programs 
on different platforms. 

\subsection{Notation and facts to remember throughout the entire text}
The following facts and notation are ``global'' and will be used freely throughout the text 
without any further references after their first occurence. 
Everything else is ``local'' to each particular (sub)section and can be safely forgotten 
when exiting the corresponding (sub)section.

\begin{itemize}
\item The definition of the Haar functions $h\ci J$ (see Level~\ref{level1});
\item The definition of the square function $Sf$ (see Level~\ref{level1});
\item The definition and the properties of the non-linear mean $M$ (see Section~\ref{2-1-1});
\item The definition and the properties of the dyadic suspension bridge $\A$ 
(see Section~\ref{2-1-2});
\item The definition of the function $\B$ (see Level~\ref{level1});
\item The notation $X(x,\tau)=\dfrac{x+\tau}{\sqrt{1-\tau^2}}$;
\item The Bellman inequality in its standard form~\eqref{BelIneq} on 
page~\pageref{BelIneq} and the inverse function form~\eqref{BelIneqInv} on 
page~\pageref{BelIneqInv};    
\item The notation 
$\wt\B(x)=\left\{\begin{aligned}1,&\quad x\le 0;\\ \tfrac 1{1+x^2},&\quad x\ge 0\end{aligned}\right.$ 
and the fact that $\wt\B$ satisfies the Bellman inequality;
\item The definition of a supersolution and the fact that $\B$ is the least supersolution (Section~\ref{2-3});
\item The notation $\Phi(x)=\int_x^{\infty}e^{-y^2/2}\,dy$ and $\Psi=\Phi^{-1}$;
\item The differential Bellman inequality $xB'(x)+B''(x)\le 0$ and its equivalence to the concavity
of the function $B\circ \Psi$ (Section~\ref{2-4});
\item The increasing property of the ratio $\dfrac{\B(x)}{\Phi(x)}$ (Section~\ref{2-4}). 
\end{itemize}

This list is here to serve as a reminder to a reader who might otherwise occasionally get lost in this text or
who might want to read its various parts in some non-trivial order.
In addition to this list, it may be useful to keep in mind the statements 
in the titles of subsections and the summary of results in Level~1 though it is not formally necessary.

\section{Level 1: Setup and main results}
\label{level1}

The celebrated Chang--Wilson--Wolff theorem (\cite{CWW}) states that, if the square function 
of a function $f$ is uniformly bounded, then $e^{a|f|^2}$ is (locally) integrable for some 
positive $a$, which, in its turn, implies that the distribution tails $\mu\{f\ge x\}$ decay 
like $e^{-ax^2}$ where $\mu$ is the usual Lebesgue measure restricted to some interval. 
This theorem holds true for both discrete and continuous versions of the square function. 
The main aim of this article is to get \textit{sharp} bounds for the distribution tails in 
the dyadic setting. 

So, let $I=[0,1]$. Let $\D$ be the collection of all dyadic subintervals of the interval $I$. 
With each dyadic interval $J\in \D$, we will associate the corresponding Haar function $h\ci J$, 
which equals $-1$ on the left half $J_-$ of the interval $J$, equals $+1$ on its right half $J_+$, 
and equals $0$ outside the interval~$J$. 

Let now $f\colon I\to\R$ be any integrable function on $I$ such that $\int_I f=0$. 
Then $f=\sum_{J\in\D}a\ci J h\ci J$ where the coefficients $a\ci J$ can be found from the formula
$a\ci J=\mu(J)^{-1}\int_I f h\ci J$ and the series converges both in $L^1$ and almost everywhere.
The dyadic square function $Sf$ of the function $f$ is then defined by the formula
$$
Sf=\sqrt{\sum_{J\in\D}a\ci J^2\chi\ci J}
$$
where $\chi\ci J=h\ci J^2$ is the characteristic function of the dyadic interval $J$.
The quantity we want to investigate is 
$$
\B(x)=\sup\{\mu\{f\ge x\}\colon\|Sf\|\ci{L^\infty}\le 1\}\,,\qquad x\in\R\,.
$$
Here is the summary of what we know and will prove in this article about the function $\B(x)$:
\begin{itemize}
\item\rule{0pt}{23pt} $\B$ is a continuous non-increasing function on $\R$;
\item\rule{0pt}{23pt} $\B(x)=1$ for all $x\le 0$ and $\B$ is strictly decreasing on $[0,+\infty)$;
\item\rule{0pt}{23pt} $\B(x)=1-\A^{-1}(x)$ for all $x\in[0,1]$ where $\A\colon[0,\frac12]\to[0,1]$ is 
the ``dyadic suspension bridge function'' constructed in the beginning of Level~\ref{level2};
\item\rule{0pt}{23pt} If $x\in[0,1]$ and $\B(x)$ is a binary rational number (i.\,e., a number of the 
kind $\frac k{2^n}$ with some non-negative integer~$k$ and~$n$), then we can explicitly construct 
a finite linear combination $f$ of Haar functions for which $\mu\{f\ge x\}=\B(x)$;
\item\rule{0pt}{23pt} There exists a positive constant $c$ (whose exact value remains unknown to us)
such that $\B(x)=c\Phi(x)$ for all $x\ge\sqrt 3$ where $\Phi$ is the Gaussian ``error function'',
i.\,e., $\displaystyle \Phi(x)=\int_x^\infty\!\! e^{-y^2/2}\,dy$.
\end{itemize}
Shortly put, this means that we know $\B$ exactly for $x\le 1$, know it up to an absolute
constant factor for $x\ge \sqrt 3$ and do not have any clear idea about what $\B$ may be between
$1$ and $\sqrt 3$.

\section{Level 2: Definitions, auxiliary results, and ideas of the proofs}
\label{level2}

\subsection{Construction of the dyadic suspension bridge function $\A$}

\subsubsection{Nonlinear mean $M$}
\label{2-1-1}
For any two real numbers $a,b$, we define their nonlinear mean $M[a,b]$ by
$$
M[a,b]=\frac{a+b}{\sqrt{4+(a-b)^2}}\,.
$$
The nonlinear mean $M[a,b]$ has the following properties.
\begin{enumerate}[(1)]
\item $M[a,a]=a$;
\item $M[a,b]=M[b,a]$;
\item $M[a,b]\le\dfrac{a+b}2$ for all $a,b\ge 0$;
\item 
$$
\frac{\partial}{\partial a}M[a,b]=\frac{4+2b^2-2ab}{[4+(a-b)^2]^{3/2}}\,;
$$ 
When $a,b\in[0,1]$, the right hand side is strictly positive and does not exceed
$\frac68=\frac34$ (the numerator is at most $6$ and the denominator is at least $8$).
It follows immediately from here that
\item $M[a,b]$ is strictly increasing in each variable in the square $[0,1]^2$ and
$M[a,b]$ lies strictly between $M[a,a]=a$ and $M[b,b]=b$ if $a,b\in[0,1]$ and $a\ne b$;
\item $|M[a,b]-a|\le\frac 34|a-b|$ for all $a,b\in [0,1]$.
\end{enumerate}

\subsubsection{Definition of $\A$}
\label{2-1-2}
Let 
$$
D_n=\left\{\frac k{2^n}\colon k=0,1,\dots,2^{n-1}\right\},\qquad n=1,2,3,\dots\,.
$$
For any $t\in D_n\setm D_{n-1}$ with $n\ge 2$, we define $t^\pm=t\pm 2^{-n}\in D_{n-1}$.
Let $D=\bigcup_{n\ge 1}D_n$ be the set of all binary rational numbers on the interval
$[0,\frac12]$. We shall define the function $\A\colon D\to[0,1]$ as follows. Put $\A(0)=0$,
$\A(\frac12)=1$. This completely defines $\A$ on $D_1$. Assume now that we already know 
the values of $\A$ on $D_{n-1}$. For each $t\in D_n\setm D_{n-1}$, we put 
$$
\A(t)=M[\A(t^-),\A(t^+)]\,.
$$
This defines $\A$ inductively on the entire $D$. The first few steps of this construction look as follows:
\begin{figure}[!h]
\hskip-240pt\vbox{\includegraphics[scale=3]{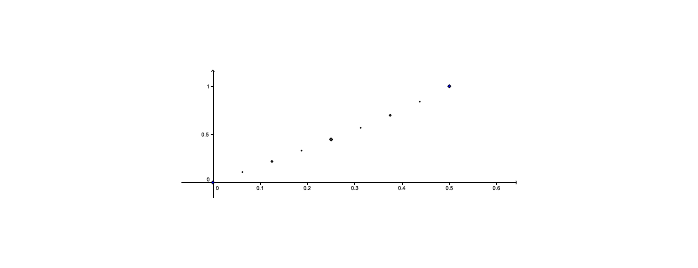}}
\caption{First steps in definition of $\A$.}
\label{fedya1}
\end{figure}

\begin{figure}[!h]
\includegraphics{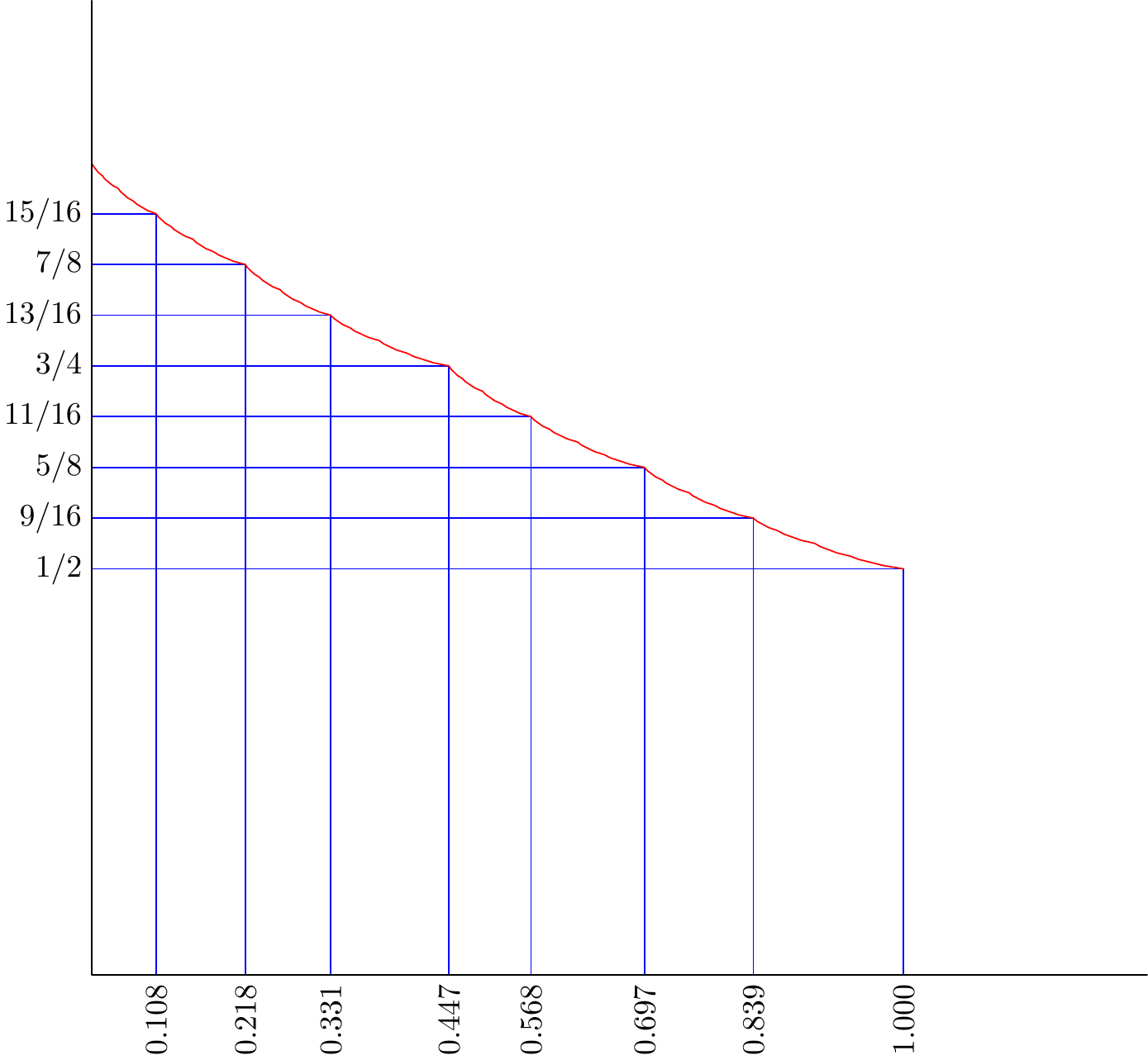}
\caption{The graph of the function $\B=1-\A^{-1}$ on $[0,1]$.}
\label{fedya1a}
\end{figure}

After completing this procedure our function $\B(x)=1-\A^{-1}(x)$ will look like it is shown
on Fig.~\ref{fedya1a}.

Property (6) of the nonlinear mean $M$ implies that the difference of values of $\A$ at any two 
neighboring points of $D_n$ does not exceed $\left(\frac 34\right)^{n-1}$. It is not hard to derive 
from here that $\A$ is uniformly continuous on $D$ and, moreover, $\A\in\ope{Lip}_\alpha$ with
$\alpha=\log_2\frac 43$. Thus, $\A$ can be extended continuously to the entire interval 
$[0,\frac 12]$. Property~(5) implies that $\A$
is strictly increasing on $D$ and, thereby, on $[0,\frac 12]$. Thus, the inverse function
$\A^{-1}\colon[0,1]\to [0,\frac 12]$ is well-defined and strictly increasing. 

\subsubsection{Properties of $\A$}
The main properties of $\A$ we shall need is the estimate
$$
\A(t)\le 2t\qquad \text{ for all }t\in[0,\tfrac12]\,,
$$
the inequality
$$
\A\left(\frac {s+t}2\right)\ge M[\A(s),\A(t)]\quad\text{ for all }s,t\in[0,\tfrac12]\,,
$$
and the fact that the function $\dfrac {\A(t)}{t}$ is non-decreasing on $(0,\frac12]$.
The first statement immediately follows from Property~(3) of the nonlinear mean $M[a,b]$ 
by induction: at the points $t=0$ and $t=1$ we have $\A(t)=2t$, and if the inequality holds
on $D_{n-1}$, then for $t\in D_n\setm D_{n-1}$ we can estimate
$$
\A(t)=M[\A(t^-),\A(t^+)]\le\frac{\A(t^-)+\A(t^+)}2\le t^-+t^+=2t\,.
$$
Thus, the assertion is true on $D$ and by continuity on the whole $[0,\tfrac12]$. 
The proofs of two other statements can be found on Level~3 in Sections~\ref{bineq} and~\ref{ratio}.

\subsection{Continuity of $\B$}
\begin{figure}[!h]
\hskip-80pt\vbox{
\includegraphics[scale=1.2]{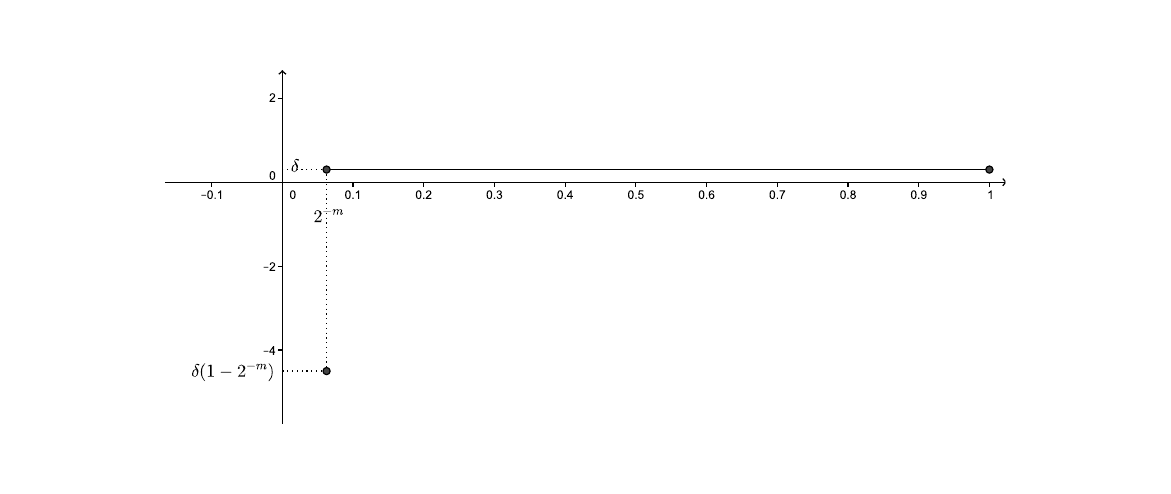}}
\caption{First summand of $g$.}
\label{fedya2}
\end{figure}

By definition, $\B$ is non-increasing on $\R$ and $\B(x)\ge 0$ for all $x\in \R$. It is easy 
to see that $\B(x)=1$ for $x\le 0$ (just consider the identically zero test-function $f$). 
Let now $x\ge 0$. Take any test-function
$f$ satisfying $\int_I f=0$ and $\|Sf\|\ci{L^\infty}\le 1$. Construct a new function
$g=g_{m,\delta}$ in the following way. Take an integer $m\ge 1$. Choose some $\delta\in (0,2^{-3m})$.
Let $I_j=[0,2^{-j}]$, $J_j=(I_j)_+=[2^{-(j+1)}, 2^{-j}]$ ($j=0,1,2,\dots$). Let $T_j$ be the 
linear mapping that maps $J_j$ onto $I$ (so, $T_0(x)=2x-1$, $T_1(x)=4x-1$, $T_2(x)=8x-1$, and so on).
Put $f_j=f\circ T_j$ on $J_j$ and $f_j=0$ on $I\setm J_j$. Now, let
$$
g(x)=\delta\sum_{j=0}^{m-1} 2^{j}h\ci{I_j}+\sqrt{1-2^{2m}\delta^2}\sum_{j=0}^m f_j\,.
$$   
The first sum may look a bit strange as written but it is just the Haar decomposition
of the function
$
\left\{
\begin{aligned}
1,\quad&\quad 2^{-m}\le x\le 1;
\\
1-2^m,&\quad 0\le x<2^{-m}
\end{aligned}
\right.
$ 
multiplied by $\delta$ (cf. Fig.~\ref{fedya2}).

Then, clearly, $\int_I g=0$. Since $f_j$ have mean $0$, are supported by disjoint dyadic intervals, and 
none of the functions $f_j$ from the second sum contains any of the function $h\ci{I_k}$ from the
first sum in its Haar decomposition, we have 
$$
(Sg)^2\le 1-2^{2m}\delta^2+\delta^2 \sum_{j=0}^{m-1} 2^{2j}\le 1
$$ 
on $I$. Finally, for each $j=0,1,\dots,m-1$, we have 
$$
\mu(\{g\ge\delta+\sqrt{1-2^{2m}\delta^2}\,x\}\cap J_j)\ge\mu(\{f_j\ge x\}\cap J_j)=2^{-(j+1)}\mu\{f\ge x\}
$$  
and, thereby, for the entire interval $I$, we have the inequality
$$
\mu\{g\ge \delta+\sqrt{1-2^{2m}\delta^2}\,x\}\ge (1-2^{-m})\mu\{f\ge x\}\,.
$$ 
Now, let us fix an integer $m>0$, thtn for every $x\in[0,2^{m-1}]$ and $\delta\in[0,2^{-3m}]$,
we have
$$
\delta+\sqrt{1-2^{2m}\delta^2}\,x\ge 
\delta+(1-2^{2m}\delta^2)x=
x+\delta(1-2^{2m}\delta x)\ge x+\frac \delta 2\,.
$$
Hence, by the definition of $\B$, we can write down the following estimate
$$
\B\left(x+\tfrac\delta 2\right)\ge 
\mu\{g\ge x+\tfrac\delta 2\}
\ge(1-2^{-m})\mu\{f\ge x\}\,.
$$
Taking the supremum over all test-functions $f$ on the right hand side, we get 
$$
\B\left(x+\tfrac\delta 2\right)\ge (1-2^{-m})\B(x)\,.
$$
Recalling that $\B$ is non-increasing and $0\le\B\le 1$, we conclude from here that
$$
0\le\B(x)-\B\left(x+\tfrac\delta 2\right)\le 2^{-m}\,,
$$
which immediately implies the uniform continuity of $\B$ on any compact subset of $\R$.

One useful corollary of this continuity result is the possibility to restrict ourselves
to the functions $f$ that are \textit{finite} linear combinations of the Haar functions in the 
definition of $\B$. Indeed, let $x\in\R$. Take any $\eps>0$. Choose $x'>x$ in such way that $\B(x')\ge
\B(x)-\eps$. Choose a function $f$ satisfying $\int_I f=0$ and $\|Sf\|\ci{L^\infty}\le 1$ such that
$\mu\{f\ge x'\}\ge\B(x')-\eps$. Let $f_n$ be the partial sums of the Haar series for $f$. Clearly,
$\int_I f_n=0$ and $Sf_n\le Sf$ everywhere on $I$. Since $f_n$ converge to $f$ almost everywhere on $I$,
we can choose $n$ such that $\mu\{f_n\ge x\}\ge \mu\{f\ge x'\}-\eps$. But then   
$\mu\{f_n\ge x\}\ge \B(x)-3\eps$. Moreover, considering the functions $g_n=(1-\frac1n)f_n$ 
instead of $f_n$, we
see that the supremum can be taken over finite linear combinations $f$ satisfying the strict inequality
$\|Sf\|\ci{L^\infty}< 1$.

\subsection{The Bellman inequality}
\label{2-3}

Take any $\tau\in(-1,1)$ and any two functions $f_-,f_+\colon I\to \R$ satisfying 
$\int_I f_\pm=0$ and $\|Sf_\pm\|\ci{L^\infty}\le 1$. Consider the function
$f$ defined by
$$
f(x)=\tau h\ci I+\sqrt{1-\tau^2}
\begin{cases}
f_-(2x),& 0\le x<\frac12\,;
\\
f_+(2x-1),& \frac12\le x\le 1\,.
\end{cases}
$$
It is easy to see that $\int_I f=0$. Also, we have
$$
((Sf)(x))^2=\tau^2+(1-\tau^2)
\begin{cases}
((Sf_-)(2x))^2,& 0\le x<\frac12\,;
\\
((Sf_+)(2x-1))^2,& \frac12\le x\le 1\,,
\end{cases}
$$ 
whence
$\|Sf\|\ci{L^\infty}\le 1$.
Now, it immediately follows from our definition of $f$ that, for every $x\in\R$,
$$
\mu\{f\ge x\}=\frac 12\left[\mu\left\{f_-\ge\frac{x+\tau}{\sqrt{1-\tau^2}}\right\}+
  \mu\left\{f_+\ge\frac{x-\tau}{\sqrt{1-\tau^2}}\right\}\right]\,.
$$
But, according to the definition of $\B$, the right hand side can be made as close to 
$\frac 12\left[\B\left(\frac{x+\tau}{\sqrt{1-\tau^2}}\right)+
  \B\left(\frac{x-\tau}{\sqrt{1-\tau^2}}\right)\right]$
as we wish by choosing appropriate~$f_\pm$. Since our function 
$f$ belongs to the class of functions over which the supremum in the 
definition of $\B(x)$ is taken, we conclude that
\begin{equation}
\label{BelIneq}
\B(x)\ge \frac 12\left[\B\left(\frac{x+\tau}{\sqrt{1-\tau^2}}\right)+
\B\left(\frac{x-\tau}{\sqrt{1-\tau^2}}\right)\right]\,.
\end{equation}
From now on, we shall use the notation $X(x,\tau)$ for $\frac{x+\tau}{\sqrt{1-\tau^2}}$. 
The inequality~\eqref{BelIneq} will be referred to as \textit{the Bellman inequality} from now on. 

We shall call every non-increasing non-negative continuous function
$B$ satisfying the Bellman inequality and the condition $B(x)=1$ for $x\le 0$ a supersolution.
Our next claim is that $\B$ is just the \textit{least} supersolution. Since $\B$ is a supersolution, 
it suffices to show that $\B(x)\le B(x)$ for any other supersolution $B$. 
It suffices to show that for any finite linear combination $f$ of the
Haar functions satisfying $\|Sf\|\ci{L^\infty}< 1$, we have $\mu\{f\ge x\}\le B(x)$ for all $x\in\R$. 
We shall prove this statement by induction on the highest level of the Haar functions in the decomposition 
of $f$ (the level of the Haar function $h\ci J$ is just the number $n$ such that $\mu(J)=2^{-n}$).
If $f$ is identically $0$ then the desired inequality immediately follows from
the definition of a supersolution. Assume that our inequality is proved for all 
linear combinations containing only Haar functions up to level $n-1$ and 
that $f$ contains only Haar functions up to level $n$.
Let $\tau$ be the coefficient at $h\ci I$ in the decomposition of $f$. Note that we must have $|\tau|<1$
(otherwise $Sf\ge 1$ on $I$). Let $T_{\pm}$ be the linear mappings that map $I$ onto $I_\pm$. 
Put $f_{\pm}=(f\circ T_\pm\mp \tau)/\sqrt{1-\tau^2}$. The functions $f_\pm$ are also finite linear
combinations of Haar functions but they contain only Haar functions up to level $n-1$
(if $n=0$, it means that $f_\pm$ are identically $0$). Also, it is not
hard to check that $\|Sf_\pm\|\ci{L^\infty}<1$. Now, clearly, 
\begin{multline*}
\mu\{f\ge x\}\le \frac 12\big[\mu\{f_-\ge X(x,-\tau)\}+\mu\{f_+\ge X(x,\tau)\}\big]
\\
\le 
\frac 12\big[B(X(x,-\tau))+B(X(x,\tau))\big]\le B(x) 
\end{multline*}
by the induction assumption and the Bellman inequality. We are done.

Now we shall characterize all triples $(x_-,x,x_+)$ of real 
numbers such that $x_\pm=X(x,\pm\tau)$ for some $\tau\in(-1,1)$. A straightforward computation
shows that in such case we must have $x=M[x_-,x_+]$ and, conversely, if $x=M[x_-,x_+]$, we 
can take $\tau=\dfrac{x_+ -x_-}{\sqrt{4+(x_+-x_-)^2}}$ and check that $x_\pm=X(x,\pm\tau)$ for this 
particular $\tau$. Thus, the Bellman inequality can be restated in the form that one must have
$$
B(x)\ge\frac12[B(x_-)+B(x_+)]
$$ 
for all triples $x_-,x,x_+$ satisfying the relation $x=M[x_-,x_+]$.

In conclusion of this section, we show that it suffices to check the Bellman inequality only 
in the case when all three numbers $x_-,x,x_+$ are non-negative. Indeed, if $x\le 0$, then
$B(x)=\max_{\R}B$ for any non-increasing function $B$ such that $B(x)=1$ for all $x\le 0$, 
and the Bellman inequality becomes trivial. If $x>0$ and, say, $x_-<0$ (note that the roles 
of $x_-$ and $x_+$ are completely symmetric), we must have $x_-=X(x,-\tau)$ with $\tau>x>0$. 
But then $X(x,\tau)>0$ and the Bellman inequality becomes stronger if we replace $\tau>x$ by 
$\tau=x$. Indeed, $B(X(x,-\tau))$ and $B(x)$ will stay the same while $B(X(x,\tau))$ will not 
decrease because $B$ is non-increasing. This remark allows us to forget about the negative 
semiaxis at all and to define a supersolution as a non-negative non-increasing continuous 
function defined on $[0,+\infty)$ and satisfying the Bellman inequality there together with 
the condition $B(0)=1$.     

\subsection{Smooth supersolutions and the differential Bellman inequality}
\label{2-4}

Suppose now that a supersolution $B$ is twice continuously differentiable
on $(0,+\infty)$. Then we have the Taylor expansion
$$
B(X(x,\pm\tau))=B(x)\pm B'(x)\tau+\frac12(xB'(x)+B''(x))\tau^2+o(\tau^2)\quad\text{as }\tau\to 0\,.
$$ 
Plugging this expansion into the Bellman inequality, we see that we must have
$$
xB'(x)+B''(x)\le 0
$$
for all $x>0$.
It is not hard to solve the corresponding linear differential equation: one possible
solution is 
$$
\Phi(x)=\int_x^\infty\!\!\! e^{-y^2/2}\,dy
$$  
and the general solution is $C_1\Phi+C_2$ where $C_1,C_2$ are arbitrary constants.

Let $\Psi\colon(0,\Phi(-\infty))\to (-\infty,+\infty)$ be the inverse function to $\Phi$. By the
inverse function theorem, we have 
$$
\Psi'=\frac 1{\Phi'\circ \Psi}=-e^{\Psi^2/2}\,.
$$  
Hence,
$$
(B\circ\Psi)''=e^{\Psi^2}((B'\circ\Psi)\cdot\Psi+ B''\circ\Psi)\,.
$$
Therefore, the differential Bellman inequality is equivalent to
concavity of $B\circ\Psi$ on $(0,\Psi(-\infty))$. Since for any non-negative
concave function $G$ on $(0,\Phi(-\infty))$, the ratio $G(t)/t$ is non-increasing,
we conclude that the ratio $\dfrac{B(\Psi(t))}{t}$ is non-increasing and, thereby,
the ratio $\dfrac {B(x)}{\Phi(x)}$ is non-decreasing on $(-\infty,+\infty)$.

The last two conditions (the concavity of $B\circ\Psi$ and the non-decreasing property
of the ratio $\dfrac {B(x)}{\Phi(x)}$) would make perfect sense for all supersolutions,
whether smooth or not. So, it would be nice to show that every supersolution
can be approximated by a $C^2$-smooth one with arbitrary precision. To do it, just
note that for every $x_-,x_+\in\R$ and every $y\ge 0$, we have 
$$
M[x_--y,x_+-y]=M[x_-,x_+]-\frac{2y}{\sqrt{4+(x_+-x_-)^2}}\ge M[x_-,x_+]-y\,.
$$ 
This allows us to conclude that if $B$ is a supersolution, then so is $B(\,\cdot\,-y)$
for all $y\ge 0$. Also note that any convex combination of supersolutions is
a supersolution as well.
Now just take any non-negative $C^2$ function $\eta$ supported by $[0,1]$ with 
total integral $1$, for $\delta>0$, define 
$\eta_\delta(x)=\delta^{-1} \eta(\delta^{-1}x)$, and consider the convolutions
$B_\delta=B*\eta_\delta$. On one hand, each $B_\delta$ is a supersolution. On the
other hand, $B_\delta\to B$ pointwise as $\delta\to\infty$.  

\subsection{$\B$ is strictly decreasing on $[0,+\infty)$}
Let us start with showing that $B(x)<1$ for all $x>0$. For this, it suffices
to note that the inequality $\|Sf\|\ci{L^\infty}\le 1$ implies 
$$
\int_I f^2=\int_I (Sf)^2\le 1\,.
$$ 
Now, if we consider the problem of maximizing $\mu\{f\ge x\}$ under the restrictions
$\int_I f=0$ and $\int_I f^2\le 1$, we shall get another function $\wt\B(x)$
on $[0,+\infty)$. Since we relaxed our restrictions, we must have $\B\le\wt\B$
everywhere. But, unlike our original
problem of finding $\B$, to find $\wt\B$ exactly is a piece of cake: we have
$$
\wt\B(x)=\frac 1{1+x^2}\qquad\text{ for all }x\ge 0\,.
$$ 
The reader can try to prove this statement himself or to look up the proof
on Level~3. Right away, we shall only mention that $\wt\B(t)$ satisfies the condition
$\wt\B(0)=1$ and the same Bellman inequality (the derivation of which is almost exactly 
the same as before; actually, the only result in this section that is impossible to repeat 
for $\wt\B$ in place of $\B$ is to show that it is the \textit{least} supersolution).  

Now, when we know that $\B(x)\le \frac1{1+x^2}<1$ for $x>0$, the strict monotonicity 
becomes relatively easy. Indeed, assume that $\B(x)=\B(y)=a$ for some $0<x<y$.
Then $a<1$. Due to the continuity of $\B$, we can choose the least $x\ge 0$ satisfying
$\B(x)=a$. This $x\ne 0$ because $\B(0)=1>a$, so we must have $x>0$. Also, we still 
have $x<y$. Take now $\tau>0$ so small that $X(x,-\tau)<x$ and $X(x,\tau)<y$. Then the Bellman
inequality immediately implies that $\B(X(x,-\tau))\le 2\B(x)-\B(X(x,\tau))\le2\B(x)-\B(y)=a$.
Since we must also have $\B(X(x,-\tau))\ge\B(x)=a$, we obtain $\B(X(x,-\tau))=a$, which contradicts
the minimality of $x$. It is worth mentioning that a 
similar argument can be used to derive continuity directly from the 
Bellman inequality. We leave the details to the reader.

The strict monotonicity property implies that $\B^{-1}$ is well defined. Also, since 
$\B(x)\le \wt\B(x)$, we must have $\B(x)\to 0+$ as $x\to\infty$. Thus, $\B^{-1}$ 
continuously maps the interval $(0,1]$ onto $[0,+\infty)$. The Bellman inequality
is equivalent to the statement that 
$$
x=\B^{-1}(\B(x))\le \B^{-1}\left(\frac{\B(x_-)+\B(x_+)}{2}\right)
$$
for all triples $x_-,\,x,\,x_+$ of non-negative numbers such that $x=M[x_-,x_+]$.
Denoting $\B(x_-)=s$, $\B(x_+)=t$, we see that the last inequality is equivalent
to 
\begin{equation*}
\label{BelIneqInv}
\B^{-1}\left(\frac {s+t}2\right)\ge M[\B^{-1}(s),\B^{-1}(t)]\,.\eqno{**}
\end{equation*}

\subsection{$\B=c\Phi$ beyond $\sqrt3$}
Our first task here will be to show that the function $\Phi$ satisfies the Bellman inequality~\eqref{BelIneq}
if $x\ge \sqrt 3$. Note that the inequality is an identity when $\tau=0$. So it suffices to show that
$$
\frac\partial{\partial \tau}[\Phi(X(x,-\tau)+\Phi(X(x,\tau))]\le 0 \qquad\text{ for all }\tau\in[0,1)\,,
$$
which, after a few simple algebraic manipulations, reduces to the inequality
$$
(1+x\tau)e^{-x\tau/(1-\tau^2)}\ge (1-x\tau)e^{x\tau/(1-\tau^2)}\,.
$$
If $x\tau\ge 1$, the left hand side is non-negative and the right hand side is non-positive.
If $x\tau<1$, we can rewrite the inequality to prove in the form 
$$
\frac 12\log\frac{1+x\tau}{1-x\tau}-\frac {x\tau}{1-\tau^2}\ge 0\,.
$$
Expanding the left hand side into a Taylor series with respect to $\tau$, we obtain the inequality
$$
\sum_{k\ge 0} x \left(\frac {x^{2k}}{2k+1}-1\right)\tau^{2k+1}\ge 0
$$
to prove. 
Observe that the coefficient at $\tau$ is always $0$ and the coefficient at $\tau^3$ is negative if
$0\le x<\sqrt 3$. It means that our inequality holds with the opposite sign 
for all sufficiently small $\tau$ if $0\le x<\sqrt 3$ and, thereby, the Bellman inequality 
fails for such $x$ and $\tau$ as well. On the other hand, if $x\ge \sqrt 3$, then \textit{all} the 
coefficients on the left hand side are non-negative and the inequality holds. 

Now let $c=\B(\sqrt 3)/\Phi(\sqrt 3)$. Consider the function $B(x)$ defined by
$$
B(x)=\left\{
\begin{aligned}
\B(x) ,&\quad x\le\sqrt 3\,;
\\
c\Phi(x),&\quad x\ge\sqrt 3\,.
\end{aligned}
\right.
$$
Note that, since the ratio $\dfrac{\B(x)}{\Phi(x)}$ is non-decreasing, we actually have
$B(x)=\min\{\B(x),c\Phi(x)\}$ everywhere on $\R$. 
Indeed, $\B(\sqrt 3)=c\Phi(\sqrt 3)$ by our choice of $c$, whence $\B\ge c\Phi$ on $[\sqrt 3,+\infty)$
and $\B\le c\Phi$ on $(-\infty,\sqrt 3]$.
Clearly, $B(x)=\B(x)=1$ for $x\le 0$,
$B$ is non-negative, continuous, and non-increasing. Let us check the Bellman inequality
for $B$. Take any triple $x_-,\,x,\,x_+$ with $x=M[x_-,x_+]$. If $x\le\sqrt 3$, we 
have 
$$
B(x)=\B(x)\ge\frac12[\B(x_-)+\B(x_+)]\ge \frac12[B(x_-)+B(x_+)]\,.
$$
If $x\ge\sqrt 3$, we 
have 
$$
B(x)=c\Phi (x)\ge\frac12[c\Phi(x_-)+c\Phi(x_+)]\ge \frac12[B(x_-)+B(x_+)]\,.
$$
Thus, $B$ is a supersolution and, therefore, $\B\le B$ everywhere. But we also know
that $\B\ge B$ everywhere. Thus, $\B=B$, i.\,e., $\B=c\Phi$ on $[\sqrt 3,\infty)$.

\subsection{$\B=1-\A^{-1}$ on $[0,1]$}

The first observation to make here is that we know the value $\B(1)$ exactly:
$\B(1)=\frac 12$. Indeed, the inequality $\B(1)\le \frac 12$ follows from the 
estimate $\B(x)\le\wt\B(x)=\frac1{1+x^2}$ and the inequality $\B(1)\ge \frac 12$
follows from the consideration of the test-function $f=h\ci I$.
Consider now the function $G(t)=\B^{-1}(1-t)$. It is continuous, increasing
and maps $[0,\frac 12]$ onto $[0,1]$. According to the
Bellman inequality in the form~\eqref{BelIneqInv}, we must have
$$
G\left(\frac{s+t}2\right)\ge M[G(s),G(t)]\qquad\text{ for all }s,t\in[0,\tfrac 12],.
$$
Also $G(0)=0=\A(0)$ and $G(\frac 12)=1=\A(\frac 12)$.
Since $M$ is monotone in each variable on $[0,1]^2$, we can easily prove by induction
that $G\ge \A$ on $D$ and, therefore, by continuity, on $[0,\frac 12]$. Applying $\B$ to both sides 
of this inequality, we conclude that $1-t\le \B(\A(t))$ on $[0,\frac 12]$. Taking $t=\A^{-1}(x)$ ($x\in[0,1]$),
we, finally, get 
$$
\B(x)\ge 1-\A^{-1}(x)\qquad\text{ for all }x\in[0,1]\,.
$$
It remains only to prove the reverse inequality. To this end, it would suffice to show
that the function
$$
B(x)=
\begin{cases}
1,& x\le 0;
\\
1-\A^{-1}(x) ,& 0\le x\le 1\,;
\\
\frac 1{1+x^2},& x\ge1\,.
\end{cases}
$$
is a supersolution. The only non-trivial property to check is the Bellman inequality. 
It has been already mentioned above that we may restrict ourselves to the case when
all three numbers $x_-,x,x_+$ are non-negative. Consider all possible cases:

\subsubsection{Case 1: all three numbers are on $[0,1]$}
In this case, we can just check the Bellman inequality in the form~\eqref{BelIneqInv}, which reduces
to the already mentioned inequality 
$$
\A\left(\frac {s+t}2\right)\ge M[\A(s),\A(t)]\quad\text{ for all }s,t\in[0,\frac12]
$$
whose proof can be found on Level~3.
\subsubsection{Case 2: $x>1$}
Here all we need is to note that, since $\A(t)\le 2t$, we have
$$
1-\A^{-1}(x)\le 1-\frac x2\le \frac 1{1+x^2}=\wt\B(x)
$$
on $[0,1]$. Therefore, we can use the fact that the Bellman inequality is true for $\wt\B$ and write
$$
B(x)=\wt \B(x)\ge\frac12[\wt\B(x_-)+\wt\B(x_+)]\ge \frac12[B(x_-)+B(x_+)]\,.
$$ 
\subsubsection{Case 3: $0<x<1$, $x_+\ge 1$}
We can always assume that it is $x_+$ that is greater than $1$ because the roles of $x_+$ and $x_-$ 
in the Bellman inequality are completely symmetric. Note that when $0<x<1$, we have 
$$
\frac\partial{\partial \tau}X(x,\tau)=\frac{1+x\tau}{(1-\tau^2)^{3/2}}>0
$$
for all $\tau\in(-1,1)$. Thus, if $x_+>x$, we must have $\tau>0$ and $x_-=X(x,-\tau)<x$.
The condition $x_-\ge 0$ implies that $\tau\le x$.

First we consider the boundary case when $x_-=0$. Then $x_+=X(x,x)=\frac{2x}{\sqrt{1-x^2}}$, 
which is greater than or equal to $1$ if and only if $x\ge \tfrac 1{\sqrt 5}$. 
Then the inequality we need to prove reduces to 
$$
B(x)\ge\frac12\bigl[B(X(x,x))+1\bigr]=\frac{1+x^2}{1+3x^2}\,.
$$
Denote that function $\frac{1+x^2}{1+3x^2}$ on the right hand side by $F(x)$ and note that 
at the endpoints of this interval we have the identities 
$B(\tfrac1{\sqrt5})=F(\tfrac1{\sqrt5})=\frac34$ and $B(1)=F(1)=\frac 12$. 
Recall also that $B\circ\Psi$ is concave on $\Psi^{-1}([\tfrac 1{\sqrt 5},1])$ 
(formally we proved this only for supersolutions but, since only arbitrarily small values
of $\tau$ were used in the proof, we can conclude that this concavity result also holds for any 
non-negative non-increasing continuous function $B$ satisfying the Bellman inequality just for the triples
$x_-,\,x,\,x_+$ contained in $[\tfrac 1{\sqrt 5},1]$). So, it would suffice to show that the function $F\circ\Psi$ 
is \textit{convex} on the same interval, which is equivalent to the assertion that $xF'(x)+F''(x)\ge 0$ 
on $[\tfrac 1{\sqrt 5},1]$. A direct computation yields
$$
xF'(x)+F''(x)=4\frac{8x^2-3x^4-1}{(1+3x^2)^3}\,.
$$  
But 
$$
8x^2-3x^4-1= 3x^2(1-x^2)+(5x^2-1)\ge 0
$$
on $[\tfrac 1{\sqrt 5},1]$ and we are done.

Now we are ready to handle the remaining case $0<x_-<x<1<x_+$. Let $\wt x_+=X(x,x)$ and
let $\wt x_-=X(x,-\tau)$ where $\tau\in(0,1)$ is chosen in such a way that $X(x,\tau)=1$.
Then $0<x_-<\wt x_-<x<1<x_+<\wt x_+$ and we have the Bellman inequality for the triples  
$0,\,x,\,\wt x_+$ and $\wt x_-,\,x,\,1$. If 
$$
B(x_+)-B(\wt x_+)\le B(0)-B(x_-) \quad\text{ or }\quad B(x_-)-B(\wt x_-)\le B(1)-B(x_+)\,,
$$  
we can prove the desired Bellman inequality for the triple $x_-,\,x,\,x_+$ by comparing 
it to the known Bellman inequality for the triple~$0,\,x,\,\wt x_+$ or~$\wt x_-,\,x,\,1$
respectively. So, the only situation that is bad for us is the one when the strict 
inequalities 
$$
B(x_+)-B(\wt x_+)>B(0)-B(x_-) \quad\text{and}\quad B(x_-)-B(\wt x_-)> B(1)-B(x_+)
$$ 
hold simultaneously. Now observe that, if four positive numbers $a,\,b,\,c,\,d$ satisfy 
$a>c$ and $b>d$, then we also have $\dfrac c{c+b}<\dfrac{a}{a+d}$. Thus, in the bad situation, we
must have
$$
\frac{B(0)-B(x_-)}{B(0)-B(\wt x_-)}<\frac{B(x_+)-B(\wt x_+)}{B(1)-B(\wt x_+)}\,.
$$

Since $\frac{\A(t)}t$ is non-decreasing on $[0,\frac12]$, we can say that
$$
\frac{B(0)-B(x_-)}{B(0)-B(\wt x_-)}\ge \frac{x_-}{\wt x_-}\,.
$$
So, in the bad situation we must have the inequality 
$$
\frac{x_-}{\wt x_-}<\frac{B(x_+)-B(\wt x_+)}{B(1)-B(\wt x_+)}\,.
$$
Note that everywhere in this inequality the function $B(x)$ coincides with $\wt\B(x)=\frac 1{1+x^2}$.
So, this is an elementary inequality (it contains fractions and square roots, of course, but still
it is a closed form inequality about functions given by explicit algebraic formulae). It turns out that
exactly the opposite inequality is always true (the proof can be found on Level~4, Subsection~\ref{4-5}), 
so we are done with this case too. 

\subsection{Optimal functions for binary rational values of $\B$}

\begin{figure}[!h]
\hskip-80pt\vbox{\includegraphics[scale=0.9]{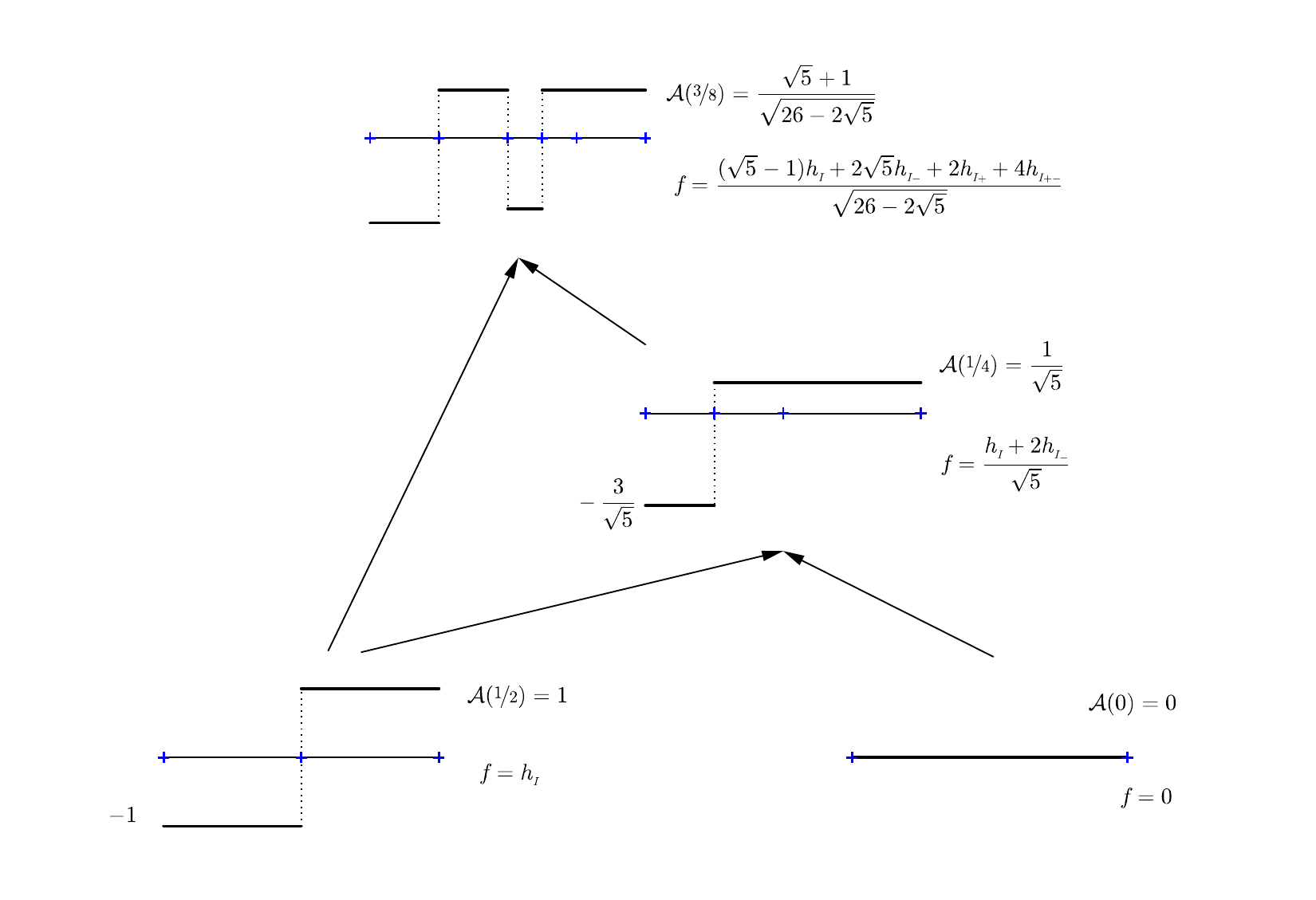}}
\caption{The construction of the extremal function for $\A(\frac 38)$.}
\label{fedya3}
\end{figure}

By the construction of the dyadic suspension bridge $\A$, for every point $t\in D\setm\{0,\frac12\}$, 
we have $\A(t)=M[\A(t_-),\A(t_+)]$. Let now $x=\A(t)$ for some $t\in D$ and let $x_-=\A(t_-)$, 
$x_+=\A(t_+)$. Then for the triple $x_-,\,x,\,x_+$, the Bellman inequality becomes an identity and 
we can say that if we have a pair $f_\pm$ of finite linear combinations of Haar functions such 
that $\|Sf_\pm\|\ci{L^\infty}\le 1$ and $\mu\{f_\mp\ge x_\pm\}=\B(x_\pm)$, then, if we take 
$\tau\in(0,1)$ such that $x_\pm=X(x,\pm\tau)$ and define $f$ by 
$$
f=\tau h\ci I+\sqrt{1-\tau^2}
\begin{cases}
f_-(2x),& 0\le x<\frac12\,;
\\
f_+(2x-1),& \frac12\le x\le 1\,,
\end{cases}
$$
we shall get a finite linear combination of Haar functions satisfying  $\|Sf\|\ci{L^\infty}\le 1$
and $\mu\{f\ge x\}=\B(x)$. Since we, indeed, have such extremal linear combinations for $x=0$ and 
$x=1$ (the identically $0$ function and the function $h\ci I$ respectively), we can now recursively 
construct an extremal linear combination for any $x=\A(t)$ with $t\in D$. Take, for instance, 
$\A(\frac 38)$. The construction of the extremal function for this value reduces to finding
the coefficient $\tau=\dfrac{\A(\frac 12)-\A(\frac 14)}{\sqrt{4+(\A(\frac 12)-\A(\frac 14))^2}}$
and two extremal functions: one for $\A(\frac 14)$ and one for $\A(\frac 12)$. The construction 
of the extremal function for $\A(\frac 14)$ reduces to finding the coefficient 
$\tau=\dfrac{\A(\frac 12)-\A(0)}{\sqrt{4+(\A(\frac 12)-\A(0))^2}}$ and two more extremal functions: 
one for $\A(0)$ and one for $\A(\frac 12)$. But we know that the extremal function for $\A(0)=0$ 
is $0$ and the extremal function for $\A(\frac 12)=1$ is $h\ci I$. So, we can put everything
together and get a linear combination of $4$ Haar functions that is extremal for $\A(\frac 38)$. 
This construction is shown on the picture~\ref{fedya3}.

The resulting linear combination is 
$$
\frac{1}{\sqrt{26-2\sqrt 5}}
\left[(\sqrt 5-1)h\ci{I}+
2\sqrt 5 h\ci {I_-}+
2h\ci{I_+}+
4h\ci{I_{+-}}
\right]\,,
$$
which, indeed, equals $\A(\frac 38)=\frac{\sqrt 5+1}{\sqrt{26-2\sqrt 5}}$ on the 
union $I_{-+}\cup I_{++}\cup I_{+-+}$ whose measure is exactly $\frac 58$.
The square function, in its turn, equals $1$ on $I_{-}\cup I_{+-}$ and is strictly
less than $1$ on $I_{++}$. 

The simplest picture is obtained when we construct an extremal function for 
$\A(2^{-n})$. What we get is just the function 
$$
\sqrt\frac{3}{4^n-1}\Big(1-2^n\chi_{_{[0,2^{-n}]}}\Big)
$$
that takes just two different values: one small positive on a big set and one large negative
on a small set. The interested reader may amuse himself with drawing more pictures, trying
to figure out how many Haar functions are needed to construct an extremal function for any
particular ``good'' value of $x$, or proving that for all other values of $x\in[0,1]$ 
there are no extremal functions at all, but we shall stop here.

\section{Level 3: Reductions to elementary inequalities} 

\subsection{$\wt\B(x)=\frac 1{1+x^2}$ for $x\ge 0$}

Recall that 
$$
\wt\B(x)\df\sup\Bigl\{\mu\{f\ge x\}\colon\int_I f=0\,,\ \int_I f^2\le 1\Bigr\}.
$$
Considering the identically zero test-function $f$, we see that $\wt\B(x)=1$ fot all $x\le 0$.
Let now $x>0$. Putting
$$
f(y)=
\begin{cases}
\displaystyle\ \ x,&\displaystyle 0\le y\le \frac 1{1+x^2};
\\
\displaystyle-\frac 1x,&\displaystyle \frac{1}{1+x^2}<y\le 1\,,
\end{cases}
$$ 
we see that $\wt\B(x)\ge \frac 1{1+x^2}$.

Now, take any test-function $f$. Let $E=\{f\ge x\}$ and let $m=\mu(E)$.
Then 
$$
\int_{I\setm E}f=-\int_E f\le -mx
$$
and 
$$
\int_{I\setm E}f^2\ge \frac1{\mu(I\setm E)}\Bigl|\int_{I\setm E}f\Bigr|^2\ge \frac{m^2x^2}{1-m}
$$
by Cauchy-Schwartz.
Thus,
$$
\int_I f^2=\int_E f^2+\int_{I\setm E}f^2\ge mx^2+\frac{m^2x^2}{1-m}=\frac{m}{1-m}x^2\,.
$$
Since this integral is bounded by $1$, we get the inequality
$$
\frac{m}{1-m}x^2\le 1\,,
$$
whence $m\le \dfrac 1{1+x^2}$\,.

One more thing we want to do in this section is to show directly that $\wt \B$ is
a supersolution. If $x, X(x,\pm\tau)\ge 0$, the Bellman inequality
$$
\wt\B(x)\ge \frac12[\wt\B(X(x,-\tau)+\wt\B(X(x,\tau))]
$$
reduces to
$$
\frac 1{1+x^2}\ge\frac12\left[\frac{1-\tau^2}{1-2x\tau+x^2}
+\frac{1-\tau^2}{1+2x\tau+x^2}\right]=\frac{(1-\tau^2)(1+x^2)}{(1+x^2)^2-4x^2\tau^2}\,,
$$
which is equivalent to
$$
(1-\tau^2)(1+x^2)^2 \le (1+x^2)^2-4x^2\tau^2\,.
$$
Subtracting $(1+x^2)^2$ from both sides, we get
$$
(1+x^2)^2\tau^2\ge 4x^2\tau^2\,.
$$
Reducing by $\tau^2$ and taking the square root of both sides,
we get the inequality
$$
1+x^2\ge 2x\,,
$$
which is obviously true.

\subsection{The inequality $\A(\tfrac{s+t}2)\ge M[A(s),A(t)]$}
\label{bineq}

Since $\A$ is continuous, it suffices to check this inequality for $s,t\in D$.
If $s,t\in D_1$, then our inequality turns into an identity. Suppose now that 
we already know that our inequality holds for all $s,t\in D_{n-1}$. To check its
validity on $D_n$, we have to consider $2$ cases:

\subsubsection{Case \textup{1:} $s\in D_n\setm D_{n-1}, t\in D_{n-1}$}
Let $s^{\pm}=s\pm 2^{-n}$. Note that $s^-$ and $s^+$ are two neighboring
points in $D_{n-1}$, whence they must lie on the same side of $t$ (it is possible that one of them
coincides with $t$). Denote $y=\A(s^-)$, $z=\A(s^+)$. By the definition of the dyadic suspension
bridge function $\A$, we then have
\begin{equation*}
\A(s)=M[y,z]\,.
\end{equation*} 
Denote $x=\A(t)$. Then 
$$
M[\A(s),\A(t)]=M[M[y,z],x]\,.
$$
Note that $\frac{s^-+t}{2}$ and  $\frac{s^++t}{2}$ are two neighboring points of $D_n$ and
the point $\frac{s+t}2\in D_{n+1}$ lies between them in the middle. Hence,
\begin{equation*}
\A\left(\frac{s+t}2\right)=
M\left[\A\left(\frac{s^-+t}2\right),\A\left(\frac{s^++t}2\right)\right]
\end{equation*} 
But, since our inequality holds on $D_{n-1}$, we have 
\begin{equation*}
\A\left(\frac{s^-+t}2\right)\ge 
M\left[\A\left(s^-\right),\A\left(t\right)\right]=M[y,x]
\end{equation*} 
and 
\begin{equation*}
\A\left(\frac{s^++t}2\right)\ge 
M\left[\A\left(s^+\right),\A\left(t\right)\right]=M[z,x]\,.
\end{equation*} 
Using monotonicity of $M$ in each argument on $[0,1]^2$, we conclude that
\begin{equation*}
\A\left(\frac{s+t}2\right)\ge 
M\left[M[z,x],M[y,x]\right]\,.
\end{equation*} 
Therefore, it would suffice to prove that
$$
M\left[M[z,x],M[y,x]\right]\ge M[M(y,z),x]
$$
for all numbers $x,y,z\in[0,1]$ such that $y$ and $z$ lie on the same side of $x$.
This will be done on Level~4 in Subsection~\ref{4-3}.

\subsubsection{Case \textup{2:} $s,t\in D_n\setm D_{n-1}$}
Without loss of generality, we may assume that $s<t$. Let, again,
$s^{\pm}=s\pm 2^{-n},\,\,t^{\pm}=t\pm 2^{-n}\in D_{n-1}$. Clearly,
$s^-<s^+\le t^-<t^+$. Denote $x=\A(s^-)$, $y=\A(s^+)$, $z=\A(t^-)$, $w=\A(t^+)$.
Then $x\le y\le z\le w$.

By the definition of the dyadic suspension
bridge function $\A$, we have
$$
\A(s)=M[x,y],\qquad \A(t)=M[z,w]\,.
$$ 
Note now that $\frac{s+t}2\in D_n$ is also a middle point for the pairs $s^-,t^+$ and $s^+,t^-$ of
the points in $D_{n-1}$. Hence, by our assumption, we have
$$
\A\left(\frac{s+t}2\right)\ge \max\{M[x,w],M[y,z]\}
$$
and, to prove the desired  inequality for $\A$ in this case, it would suffice to show
that
$$
M[M[x,y],M[z,w]]\le \max\{M[x,w],M[y,z]\}\,,
$$
provided that $0\le x\le y\le z\le w\le 1$.
This will be done on Level~4 in Subsection~\ref{4-4}.

\subsection{The ratio $\A(t)/t$ increases}
\label{ratio}

Since $\A$ is continuous, it suffices to check this property for $t\in D$.
We shall show by induction on $m$ that, for every $t_0\in D_n\setm\{\frac12\}$, the
ratio $\frac{A(t)-A(t_0)}{t-t_0}$ is non-decreasing on $D_{n+m}\cap(t_0,t_0+2^{-n}]$.
The property to prove coincides with this statement for $n=1$, $t_0=0$.

The base of induction $m=1$ is fairly simple. The interval \hbox{$(t_0,t_0+2^{-n}]$}
contains just two points of $D_{n+1}$: $t_1=t_0+2^{-(n+1)}\in D_{n+1}\setm D_n$ 
and $t_2=t^0+2^{-n}\in D_n$. By the definition of $\A$ and property $(3)$ of $M$, 
we have
$$
A(t_1)=M[A(t_0),A(t_2)]\le \frac{A(t_0)+A(t_2)}{2}\,,
$$
whence 
$$
\frac{A(t_1)-A(t_0)}{t_1-t_0}\le \frac{A(t_2)-A(t_0)}{2(t_1-t_0)}=\frac{A(t_2)-A(t_0)}{t_2-t_0}\,.
$$

Assume now that the statement is already proved for $m-1\ge 1$. Let $t_0\in D_n$ and let, again,
$t_1=t_0+2^{-(n+1)}\in D_{n+1}\setm D_n$ 
and $t_2=t^0+2^{-n}\in D_n\subset D_{n+1}$.
By the induction assumption applied to $n+1$ and $m-1$ instead of $n$ and $m$, we see that
the ratio
$\frac{A(t)-A(t_0)}{t-t_0}$ is non-decreasing on $D_{n+m}\cap(t_0,t_1]$
and the ratio
$\frac{A(t)-A(t_1)}{t-t_1}$ is non-decreasing on $D_{n+m}\cap(t_1,t_2]$.
Note also that, for $t\in (t_1,t_2]$, we have the identity
$$
\frac{A(t)-A(t_0)}{t-t_0}=
\frac{A(t_1)-A(t_0)}{t_1-t_0}+
\frac{t-t_1}{t-t_0}
\left[
\frac{A(t)-A(t_1)}{t-t_1}-
\frac{A(t_1)-A(t_0)}{t_1-t_0}
\right]\,.
$$ 
Since $t\mapsto \frac{t-t_1}{t-t_0}$ is a positive increasing function on $(t_1,t_2]$, checking
the non-decreasing property of the ratio $\frac{A(t)-A(t_0)}{t-t_0}$ reduces to showing that the factor
$\frac{A(t)-A(t_1)}{t-t_1}-
\frac{A(t_1)-A(t_0)}{t_1-t_0}$ is non-negative and non-decreasing on $D_{n+m}\cap(t_1,t_2]$. We know that 
it is non-decreasing by the induction assumption and, therefore, it suffices to check its non-negativity 
at the least element of $D_{n+m}\cap(t_1,t_2]$, which is $t'=t_1+2^{-(n+m)}$.

Let $x=\A(t_1)$. By the construction of the function $\A$, we have $\A(t')=y_m$ where the sequence
$y_j$ is defined recursively by $y_1=\A(t_2)$, $y_j=M[x,y_{j-1}]$ for all $j\ge 2$.
We shall also consider the auxiliary sequence $z_j$  
defined recursively by $z_1=\A(t_0)$, $z_j=\frac{z_{j-1}+x}2$ for all $j\ge 2$.

Note that
$$
\frac{\A(t_1)-A(t_0)}{t_1-t_0}=\frac {x-z_1}{t_1-t_0}=2^{m-1}\frac {x-z_m}{t_1-t_0}\,.
$$
Also,
$$
\frac{\A(t')-A(t_1)}{t'-t_1}=2^{m-1}\frac {y_m-x}{t_2-t_1}\,.
$$
Since $t_2-t_1=t_1-t_0=2^{-(n+1)}$, our task reduces to proving that
$y_m-x\ge x-z_m$ or, equivalently, $x\le \frac{z_m+y_m}2$. We shall show by induction on $j$ that
even the stronger inequality $x\le M[z_j,y_j]$ holds for all $j\ge 1$. 

For the base we have the identity $x=M[z_1,y_1]$ following right from the definition
of $\A$ (recall that $t_0$ and $t_2$ are two neighboring points of $D_n$ and
$t_1\in D_{n+1}$ lies in the middle between them).

To make the induction step, it would suffice to show that for every triple
$0\le z\le x\le y\le 1$ satisfying $x\le M[z,y]$, we also have 
$$
x\le M\left[\frac{x+z}2,M[x,y]\right]\,.
$$
Unfortunately, we have managed to prove it only under the additional restriction $y-z\le\frac34$.
Fortunately, this restriction holds automatically almost always. If $n\ge 2$, then using 
property~$(6)$ of the nonlinear mean, we get
$$
y_j-z_j\le y_1-z_1=A(t_2)-A(t_0)\le \left(\frac34\right)^{n-1}\le \frac34 
$$
for all $j\ge 1$.
Also, if $j\ge 2$, we have
$$
y_j-z_j=M[x,y_{j-1}]-\frac{z_{j-1}\!+\!x}{2}\le 
\frac{x\!+\!y_{j-1}}{2}-\frac{z_{j-1}\!+\!x}{2}=\frac{y_{j-1}\!-\!z_{j-1}}{2}\le \frac12
$$
for all $n\ge 1$.

Thus, the only case we cannot cover by our induction step is $n=1$, $j=2$. 
We will have to add it to the base. It is just the numerical inequality
$$
\frac1{\sqrt5}\le M\left[\frac1{2\sqrt5},  M\left[\frac1{\sqrt5},1\right]\right]\,,
$$
which shall be checked on Level~5.

The last observation we want to make in this section is that, instead of checking the inequality
$x\le M[\frac{x+z}2, M[x,y]]$ for all triples $0\le z\le x\le y\le 1$ satisfying
$y-z\le \frac 34$, $x\le M[z,y]$, we can check it only for the case $0\le y-z\le \frac34$, $x=M[z,y]$.
Indeed, since $M[z,y]\ge x$, $M[z,x]\le M[x,x]=x$, and $M$ is continuous, we can use the intermediate
value theorem and find $y'\in[x,y]$ such that $M[z,y']=x$. Obviously, $y'-z\le y-z\le\frac 34$ too.
Now, if we know that $x\le M[\frac{x+z}2, M[x,y']]$, we can just use monotonicity of $M$ twice and conclude
that $x\le M[\frac{x+z}2, M[x,y]]$ as well. This observation allows to eliminate $x$ from the inequality to
prove altogether. All we need to show is that
$$
M[z,y]\le M\left[\frac{z+M[z,y]}{2}, M[M[z,y],y]\right]
$$  
whenever $0\le z\le y\le 1$ and $y-z\le\frac 34$.
This will be done on Level~4 in Section~\ref{4-6}

\section{Level 4: Proofs of elementary inequalities}

\subsection{General idea}

We shall reduce all our elementary inequalities to checking non-negativity of some
polynomials of $2$ or $3$ variables with rational coefficients on the unit square
$[0,1]^2$ or the unit cube $[0,1]^3$. Since the polynomials that will arise on this
way are quite large (typically, they can be presented on $1$ or $2$ pages, but one of them,
if written down in full, would occupy more than $6$ pages), to check their non-negativity
by hand would be quite a tedious task, to say the very least. So, we will need some simple and
easy program to test for non-negativity that would allow us to delegate the actual work to
a computer.

\subsection{Non-negativity test}

We shall start with polynomials of one variable. Suppose that we want to check
that $P(x)=a_0+a_1 x+a_2 x^2+\dots+a_n x^n\ge 0$ on $[0,1]$. 
Then, of course, we should check, at least that $a_0=P(0)\ge 0$.
Suppose it is so. Write our polynomial in the form
$$
P(x)=a_0+x(a_1 +a_2 x+\dots+a_n x^{n-1})
$$ 
and replace the first factor $x$ by $x_1$. We shall get a polynomial
of $2$ variables 
$$
Q(x,x_1)=a_0+x_1(a_1 +a_2 x+\dots+a_n x^{n-1})\,.
$$ 
Clearly, if $Q$ is non-negative on $[0,1]^2$, then $P$ is non-negative on $[0,1]$.
But $Q$ is linear in $x_1$, so it suffices to check its non-negativity at the 
endpoints $x_1=0$ and $x_1=1$. The first case reduces to checking that $a_0\ge 0$, which has 
been done already, and the second case reduces to checking the non-negativity
of the polynomial
$$
(a_0+a_1) +a_2 x+\dots+a_n x^{n-1}=P(0)+ \frac {P(x)-P(0)}x\,,
$$
which is a polynomial of smaller degree.

This observation leads to the following informal algorithm:

\bigskip
\begin{enumerate}[(1)]
\item Is $P(0)\ge 0$? If not, stop and report failure. If yes, proceed.
\item Is $P$ constant? If yes, stop and report success. If no, proceed.
\item Replace $P$ by $P(0)+\frac {P(x)-P(0)}x$ and go back to step (1).
\end{enumerate}

\bigskip
Of course, since we know the number of steps needed to reduce
the polynomial to a constant exactly (it is just the degree of the polynomial), 
the ``go to'' operation will be actually replaced by a ``for'' loop in the real
program. Otherwise the algorithm we shall use is exactly as written. Here is
the formal program for Mathematica the reader may want to play with a bit before
proceeding just to make sure it works as promised.

\bigskip
\begin{verbatim}
P[x_]=...;
flag=False;
n=Exponent[P[x],x];
For[k=0, k<n+1, k++,
    If[P[0]<0, flag=True; Break[] ];
    P[x_]=Expand[P[0]+(P[x]-P[0])/x]
   ];
If[flag, Print["Test failed"], Print["Test successful"]];   
\end{verbatim}

\bigskip\noindent\hglue0pt
Of course, when running this program, 
instead of three dots, one needs to plug in the polynomial one wants to test.
Also, the reader may want to execute the command

\bigskip
\begin{verbatim}
Clear[P,x,n,k,flag];
\end{verbatim}

\bigskip\noindent\hglue0pt
prior to running this program in Mathematica if he has already introduced
the corresponding variables during his previous work. Note, by the way, that,
while the initialization of $P$ in the beginning can be done by the $:=$ operator
instead of $=$, using $:=$ for modifying $P$ inside the loop 
will result in an infinite recursion,
which can effectively suspend the operations of a computer.
So, when copying this and other programs of ours from the paper, one should
pay attention to various ``minor'' details like this one.

If one thinks a bit about what this test really does, one can realize that what
is actually checked is the non-negativity of the polyaffine form
$$
Q(x_1,x_2,\dots,x_n)=a_0+a_1 x_1+a_2 x_1 x_2+\dots+a_n x_1 x_2\dots x_n
$$
on $[0,1]^n$ and the test really reduces to checking that all partial sums of
the coefficients starting with $a_0$ are non-negative. In this form, the test
is well-known to any analyst in the form of the statement that non-negativity
of Ces\`aro partial sums implies non-negativity of Abel--Poisson ones. What is 
surprising here is not the test itself, but its uncanny effectiveness.

The test can easily be generalized to polynomials of more than one variable.
All we need to do is to treat a polynomial of $2$ or more variables as a polynomial
of one fixed variable with coefficients that are polynomials of other variables.
In this way, checking the non-negativity of one polynomial of, say, $3$ variables
is reduced to checking non-negativity of several polynomials of $2$ variables, to
each of which we can apply our test again. It seems that the best way to program
such a test is to write a recursive subroutine but, since the number of variables
in all our applications does not exceed $3$ and since the sleekness of our programming 
was the least of our concerns when working on this project, we just wrote the 
test for $3$ variables as follows:

\bigskip
\begin{verbatim}
LinearTest=Function[
  flag=False;
  nz=Exponent[R[x,y,z],z];
  For[kz=0, kz<nz+1, kz++,
      S[x_,y_]=R[x,y,0];
      ny=Exponent[S[x,y],y];
      For[ky=0, ky<ny+1, ky++,
          T[x_]=S[x,0];
          nx=Exponent[T[x],x];
          For[kx=0, kx<nx+1, kx++,
              If[T[0]<0, flag=True; Break[] ];
              T[x_]=Expand[T[0]+(T[x]-T[0])/x]
             ] 
          If[flag, Break[] ];   
          S[x_,y_]=Expand[S[x,0]+(S[x,y]-S[x,0])/y]
         ]; 
     If[flag, Break[] ];   
     R[x_,y_,z_]=Expand[R[x,y,0]+(R[x,y,z]-R[x,y,0])/z]
     ];
  If[flag, Print["Test failed"], Print["Test succeded"] ]; 
]
\end{verbatim}

\bigskip
The way to apply the test to some actual polynomial is to execute 
the sequence of commands 

\bigskip
\begin{verbatim}
R[x,y,z]=...;
LinearTest[];
\end{verbatim}

\bigskip\noindent\hglue0pt
where, again, three dots should be replaced by the actual polynomial one
wants to test. Note that we can interpret a polynomial of fewer than
three variables as a polynomial of three variables, so this three-variable
test can be applied verbatim to polynomials of $2$ variables as well with
the same syntax. Again, what is actually checked is the non-negativity
of a polyaffine form and the test reduces to checking that all the rectangular
partial sums of the coefficients are non-negative (the last observation implies,
in particular, that the order in which the variables are used in the test is
of no importance; we make this remark because we ourselves were stupid
enough to apply the test with all~$6$ possible rearrangements
of variables $x,y,z$ before getting convinced that it fails). On the 
other hand, it is quite possible that the test will fail for $P(x)$, but 
will succeed for $P(1-x)$: just consider $4x-6x^2+4x^3-x^4=1-(1-x)^4$. So,
some clever fiddling with variables may occasionally help.

\subsection{The inequality $M[M[y,x],M[z,x]]\ge M[M[y,z],x]$}
\label{4-3}

Recall that we need to prove this inequality under the assumptions that $x,y,z\in[0,1]$,
and $y$ and $z$ lie on the same side of $x$. Denote 
$$
a=M[y,x]\,,\qquad b=M[z,x]\,, \qquad c=M[y,z]\,.
$$
Raising both sides of the original inequality to the second power (which is 
legitimate because they are non-negative), we see that we need to show that
$$
\frac{(a+b)^2}{4+(a-b)^2}\ge \frac{(x+c)^2}{4+(x-c)^2}\,.
$$
Multiplying by the denominators, we can rewrite it as
$$
(a+b)^2\left(4+(x-c)^2\right)\ge (x+c)^2\left(4+(a-b)^2\right)\,.
$$
Opening the parentheses and regrouping the terms, we get
\begin{multline*}
\left((a^2+b^2)+2ab\right)\left((4+x^2+c^2)-2xc\right)
\\
\ge   
\left((x^2+c^2)+2xc\right)\left((4+a^2+b^2)-2ab\right)\,.            
\end{multline*}
Putting all the terms containing the product $ab$ on the left and all other terms
on the right, we get the inequality
$$
2ab(4+2x^2+2c^2)\ge 4(x^2+c^2-a^2-b^2)+2xc(4+2a^2+2b^2)\,,
$$ 
which, after division by $4$, reduces to 
$$
ab(2+x^2+c^2)\ge (x^2+c^2-a^2-b^2)+xc(2+a^2+b^2)\,.
$$ 
Now denote 
$$
U=2+x^2+c^2\,,\qquad V=x^2+c^2-a^2-b^2\,,\qquad W=x(2+a^2+b^2)\,.
$$
Our inequality becomes 
$$
Uab\ge V+Wc\,.
$$
Since the left hand side is, clearly, non-negative, it suffices to check 
the squared inequality 
$$
U^2a^2b^2\ge V^2+W^2c^2+2VWc\,,
$$ 
or, which is the same, 
$$
U^2a^2b^2-V^2-W^2c^2\ge 2VWc\,. 
$$
At this point, we need information about the sign of the left hand side
$F=U^2a^2b^2-V^2-W^2c^2$ to proceed. Note that $a^2, b^2, c^2$ are 
\textit{rational} functions of $x,y,z$ and, therefore, so are $U,V,W$ and $F$.
We can program the computation of $F$ in Mathematica as follows: 

\bigskip
\begin{verbatim}
Den[x_,y_]=4+(x-y)^2;
MM[x_,y_]=(x+y)^2/Den[x,y];
U[x_,y_,z_]=2+x^2+MM[y,z];
V[x_,y_,z_]=x^2+MM[y,z]-MM[y,x]-MM[z,x];
W[x_,y_,z_]=(2+MM[y,x]+MM[z,x])*x;
F[x_,y_,z_]=U[x,y,z]^2*MM[y,x]*MM[z,x]-V[x,y,z]^2-W[x,y,z]^2*MM[y,z];
Print[Factor[F[x,y,z]]];
\end{verbatim}

\bigskip\noindent\hglue0pt
The output looks like
$$
\begin{aligned}
&(4(x - y)(x - z)(1024x y + 1536 x^3 y +\dots 
+2 y^3 z^7- 2 x^2 y^3 z^7))/
\\
&((4 + x^2 - 2 x y + y^2)^2
(4 + x^2 - 2 x z + z^2)^2(4 + y^2 - 2yz + z^2)^2)\,.
\end{aligned}
$$
Since $(x-y)(x-z)\ge 0$ under our assumptions and the product
of the denominators is obviously positive, we only need to determine
the sign of the huge polynomial $P_1(x,y,z)$ in the middle (if written in full, it
occupies about half-page:

\medskip\noindent\hglue0pt
{\tiny
$P_1(x,y,z)=1024 x y + 1536 x^3 y + 832 x^5 y + 192 x^7 y + 16 x^9 y - 
    512 x^2 y^2 - 624 x^4 y^2 - 248 x^6 y^2 - 32 x^8 y^2 + 
    1024 x y^3 + 1344 x^3 y^3 + 624 x^5 y^3 + 120 x^7 y^3 + 
    8 x^9 y^3 - 316 x^2 y^4 - 334 x^4 y^4 - 114 x^6 y^4 - 
    12 x^8 y^4 + 316 x y^5 + 324 x^3 y^5 + 117 x^5 y^5 + 
    18 x^7 y^5 + x^9 y^5 - 45 x^2 y^6 - 45 x^4 y^6 - 
    13 x^6 y^6 - x^8 y^6 + 27 x y^7 + 18 x^3 y^7 + 
    3 x^5 y^7 + 1024 x z + 1536 x^3 z + 832 x^5 z + 
    192 x^7 z + 16 x^9 z - 1024 x^2 y z - 1312 x^4 y z - 
    528 x^6 y z - 64 x^8 y z + 512 x y^2 z + 576 x^3 y^2 z + 
    208 x^5 y^2 z + 8 x^7 y^2 z - 8 x^9 y^2 z - 
    624 x^2 y^3 z - 544 x^4 y^3 z - 112 x^6 y^3 z + 
    12 x y^4 z - 44 x^3 y^4 z - 39 x^5 y^4 z - 
    14 x^7 y^4 z - 3 x^9 y^4 z + 20 y^5 z - 150 x^2 y^5 z - 
    92 x^4 y^5 z - 4 x^6 y^5 z + 2 x^8 y^5 z - 
    11 x y^6 z + 9 x^5 y^6 z + 2 x^7 y^6 z + 9 y^7 z - 
    21 x^2 y^7 z - 11 x^4 y^7 z - x^6 y^7 z - 512 x^2 z^2 - 
    624 x^4 z^2 - 248 x^6 z^2 - 32 x^8 z^2 + 512 x y z^2 + 
    576 x^3 y z^2 + 208 x^5 y z^2 + 8 x^7 y z^2 - 
    8 x^9 y z^2 - 680 x^2 y^2 z^2 - 292 x^4 y^2 z^2 + 
    68 x^6 y^2 z^2 + 24 x^8 y^2 z^2 + 504 x y^3 z^2 + 
    552 x^3 y^3 z^2 + 146 x^5 y^3 z^2 - 4 x^7 y^3 z^2 + 
    2 x^9 y^3 z^2 - 251 x^2 y^4 z^2 - 35 x^4 y^4 z^2 + 
    21 x^6 y^4 z^2 + x^8 y^4 z^2 + 111 x y^5 z^2 + 
    92 x^3 y^5 z^2 + 11 x^5 y^5 z^2 - 6 x^7 y^5 z^2 - 
    4 y^6 z^2 - 24 x^2 y^6 z^2 + 2 x^4 y^6 z^2 + 
    2 x^6 y^6 z^2 + x y^7 z^2 + 6 x^3 y^7 z^2 + 
    x^5 y^7 z^2 + 1024 x z^3 + 1344 x^3 z^3 + 624 x^5 z^3 + 
    120 x^7 z^3 + 8 x^9 z^3 - 624 x^2 y z^3 - 
    544 x^4 y z^3 - 112 x^6 y z^3 + 504 x y^2 z^3 + 
    552 x^3 y^2 z^3 + 146 x^5 y^2 z^3 - 4 x^7 y^2 z^3 + 
    2 x^9 y^2 z^3 - 40 y^3 z^3 - 132 x^2 y^3 z^3 - 
    40 x^4 y^3 z^3 - 8 x^6 y^3 z^3 - 4 x^8 y^3 z^3 + 
    65 x y^4 z^3 + 18 x^3 y^4 z^3 - 23 x^5 y^4 z^3 + 
    4 x^7 y^4 z^3 - 9 y^5 z^3 + 21 x^2 y^5 z^3 + 
    11 x^4 y^5 z^3 + x^6 y^5 z^3 + 5 x y^6 z^3 + 
    6 x^3 y^6 z^3 - 3 x^5 y^6 z^3 + 2 y^7 z^3 - 
    2 x^2 y^7 z^3 - 316 x^2 z^4 - 334 x^4 z^4 - 
    114 x^6 z^4 - 12 x^8 z^4 + 12 x y z^4 - 44 x^3 y z^4 - 
    39 x^5 y z^4 - 14 x^7 y z^4 - 3 x^9 y z^4 - 
    251 x^2 y^2 z^4 - 35 x^4 y^2 z^4 + 21 x^6 y^2 z^4 + 
    x^8 y^2 z^4 + 65 x y^3 z^4 + 18 x^3 y^3 z^4 - 
    23 x^5 y^3 z^4 + 4 x^7 y^3 z^4 + 8 y^4 z^4 - 
    80 x^2 y^4 z^4 - 4 x^4 y^4 z^4 - 4 x^6 y^4 z^4 + 
    10 x y^5 z^4 - 12 x^3 y^5 z^4 + 2 x^5 y^5 z^4 - 
    8 x^2 y^6 z^4 + 316 x z^5 + 324 x^3 z^5 + 117 x^5 z^5 + 
    18 x^7 z^5 + x^9 z^5 + 20 y z^5 - 150 x^2 y z^5 - 
    92 x^4 y z^5 - 4 x^6 y z^5 + 2 x^8 y z^5 + 
    111 x y^2 z^5 + 92 x^3 y^2 z^5 + 11 x^5 y^2 z^5 - 
    6 x^7 y^2 z^5 - 9 y^3 z^5 + 21 x^2 y^3 z^5 + 
    11 x^4 y^3 z^5 + x^6 y^3 z^5 + 10 x y^4 z^5 - 
    12 x^3 y^4 z^5 + 2 x^5 y^4 z^5 - 4 y^5 z^5 + 
    20 x^2 y^5 z^5 - 45 x^2 z^6 - 45 x^4 z^6 - 13 x^6 z^6 - 
    x^8 z^6 - 11 x y z^6 + 9 x^5 y z^6 + 2 x^7 y z^6 - 
    4 y^2 z^6 - 24 x^2 y^2 z^6 + 2 x^4 y^2 z^6 + 
    2 x^6 y^2 z^6 + 5 x y^3 z^6 + 6 x^3 y^3 z^6 - 
    3 x^5 y^3 z^6 - 8 x^2 y^4 z^6 + 27 x z^7 + 
    18 x^3 z^7 + 3 x^5 z^7 + 9 y z^7 - 21 x^2 y z^7 - 
    11 x^4 y z^7 - x^6 y z^7 + x y^2 z^7 + 6 x^3 y^2 z^7 + 
    x^5 y^2 z^7 + 2 y^3 z^7 - 2 x^2 y^3 z^7 $}

\medskip
To recover $P_1$ from $F$, it is enough to execute the command

\bigskip
\begin{verbatim}
P1[x_,y_,z_]=Factor[F[x,y,z]*
         Den[x,y]^2*Den[x,z]^2*Den[y,z]^2/4/(x-y)/(x-z)];
\end{verbatim}

\bigskip\noindent\hglue0pt
If we apply our non-negativity test to the polynomial 
$P_1(x,y,z)$ directly, then it reports failure. But after we looked
into how exactly it failed, we discovered that it fails already on the 
polynomial $P(0,y,z)$. This particular polynomial is not hard to factor:
executing the command 

\bigskip
\begin{verbatim}
Print[Factor[P1[0,y,z]]];
\end{verbatim}

\bigskip\noindent\hglue0pt
we get
$$
P_1(0,y,z)=y(y-z)^2z(y+z)^2(20+9y^2-4yz+9z^2+2y^2z^2)\,,
$$
which is obviously a non-negative function on $[0,1]^2$.
So, it will suffice to show that $P_1(x,y,z)-P_1(0,y,z)$ is non-negative
and that can be done by our test: the execution of the commands

\bigskip
\begin{verbatim}
R[x_,y_,z_]=P1[x,y,z]-P1[0,y,z];
LinearTest[];
\end{verbatim}

\bigskip\noindent\hglue0pt
reports a successful completion of the test.

Now, once we know that $F\ge 0$, we can say that our inequality would follow from
the squared inequality
$$
F^2-4 V^2 W^2 c^2\ge 0
$$
whose left hand side is a rational function of $x,y,z$. Remembering
that we had trouble with $x=0$ last time, we should expect it again
because $W=(2+a^2+b^2)x$ has a factor $x$ in it, which means that our rational 
function and the corresponding huge polynomial factor in it are the same as in $F^2$
when $x=0$. Fortunately, this time we do not even need to factor anything to 
realize that $F^2\ge 0$ when $x=0$: the square is always non-negative. Let us keep
it in mind and execute the commands

\bigskip
\begin{verbatim}
G[x_,y_,z_]=F[x,y,z]^2-4*V[x,y,z]^2*W[x,y,z]^2*MM[y,z];
Print[Factor[G[x,y,z]]];
\end{verbatim}

\bigskip\noindent\hglue0pt
The output looks like 
$$
\begin{aligned}
&(16(x - y)^2(x - z)^2(y-z)^2
\\
&(32768x^5y + 65536 x^7 y +\dots 
+12 x^4 y^6 z^{12}- 4 x^6 y^6 z^{12}))/
\\
&((4 + x^2 - 2 x y + y^2)^4
(4 + x^2 - 2 x z + z^2)^4(4 + y^2 - 2yz + z^2)^4)
\end{aligned}
$$
with a huge polynomial $P_2(x,y,z)$ in the middle (about $6$
times as long as $P_1$). Now we know that 
$P_2(0,y,z)=P_1(0,y,z)^2\ge 0$ and that we may have some trouble at this
level. So, we will immediately subtract $P_2(0,y,z)$ and apply our non-negativity
test to the difference. The corresponding sequence of commands to execute is the 
following:

\bigskip 
\begin{verbatim}
P2[x_,y_,z_]=Factor[G[x,y,z]*
             Den[x,y]^4*Den[x,z]^4*Den[y,z]^4/
                        16/(x-y)^2/(x-z)^2/(y-z)^2];
R[x_,y_,z_]=P2[x,y,z]-P2[0,y,z];
LinearTest[];
\end{verbatim}

\bigskip\noindent\hglue0pt
The test reports success, thus finishing the proof.

\subsection{The inequality $M[M[x,y],M[z,w]]\le \max\{M[x,w],M[y,z]\}$}
\label{4-4}

Denote the right hand side by $u$. Since $x\le y\le u\le z\le w$, we can raise
$x$ to $x'\in[x,u]$ and $y$ to $y'\in[y,u]$ such that
$M[x',w]=M[y',z]=u$. The right hand side will not change and the left hand side will not decrease,
so the inequality will get only stronger.

Now, choose $\sigma$ and $\tau$ such that
$$
x'=X(u,-\sigma)\,,\quad
w=X(u,\sigma)\,,\qquad
y'=X(u,-\tau)\,,\quad
z=X(u,\tau)\,.
$$
Since $u\le z\le w\le 1$, we must have $0\le\tau\le\sigma$ (recall that
$\frac\partial{\partial\tau}X(u,\tau)=\frac{1+u\tau}{(1-\tau^2)^{3/2}}>0$
when $0\le u\le 1$, $|\tau|<1$). Our inequality can be rewritten as
$$
M\big[M[X(u,-\sigma), X(u,-\tau) ], M[X(u,\sigma), X(u,\tau) ]  \big]\le u\,.
$$
Since the expressions $X(u,\pm\sigma)=\frac{u\pm\sigma}{\sqrt{1-\sigma^2}}$ and 
$X(u,\pm\tau)=\frac{u\pm\tau}{\sqrt{1-\tau^2}}$ contain square roots and we would 
strongly prefer to deal with purely rational functions, we will make one more change
of variable and put $\sigma=\frac{2s}{1+s^2}$, $\tau=\frac{2t}{1+t^2}$ ($s,t\in[0,1)$).
Then 
$$
X(u,\pm\sigma)=\frac{(1+s^2)u\pm 2s}{1-s^2}\quad\text{ and }\quad X(u,\pm\tau)=\frac{(1+s^2)u\pm 2t}{1-t^2}\,.
$$
Now it is time to discuss the possible joint range of the variables $u,s,t$. Since the 
function $r\mapsto \frac{2r}{1+r^2}$ is strictly increasing on $[0,1)$, we must have
$0\le t\le s$ because $0\le\tau\le\sigma$. Also, $\sigma=\frac{w-x'}{\sqrt{4+(w-x')^2}}$
and, since $0\le w-x'\le 1$ and since the function $r\mapsto \frac r{\sqrt{4+r^2}}$ is 
increasing on $[0,+\infty)$, we get $\sigma\le\frac{1}{\sqrt 5}<\frac 8{17}$, whence
$s\le\frac 14$. 

As to $u$, since $x'\ge 0$, we have $u\ge \sigma=\frac{2s}{1+s^2}$ and, surely, $u\le 1$.
Thus, the joint range of our variables is contained in the domain
$$
0\le t\le s\le \frac14\,\qquad \frac{2s}{1+s^2}\le u\le 1\,.
$$
Now we are ready to proceed with the proof. Denote 
$$
M_\pm=M[X(u,\pm\sigma),X(u,\pm\tau)]
$$ 
and square both sides of the inequality. We get
$$
\frac{(M_-+M_+)^2}{4+(M_--M_+)^2}\le u^2\,,
$$
which can be rewritten as 
$$
(M_-+M_+)^2\le \left(4+(M_--M_+)^2\right)u^2\,,
$$
or, after opening the parentheses and regrouping the terms, as
$$
2M_-M_+(1+u^2)\le 4u^2-(1-u^2)(M_-^2+M_+^2)\,.
$$
Observe now that 
$$
M_\pm^2\le \frac{(X(u,\pm\tau)+X(u,\pm\sigma))^2}{4}\,, 
$$
and 
$$
X(u,\pm\tau)=\frac{u\pm\tau}{\sqrt{1-\tau^2}}\le \frac{u\pm\tau}{\sqrt{1-u^2}}
$$
and a similar inequality holds for $X(u,\pm\sigma)$. Thus
\begin{multline*}
4u^2-(1-u^2)(M_-^2+M_+^2)\ge 4u^2-\frac{(2u-\sigma-\tau)^2}{4}
-\frac{(2u+\sigma+\tau)^2}{4}
\\
=4u^2-2u^2-\frac{(\sigma+\tau)^2}{2}
\ge 4u^2-2u^2-\frac{(2u)^2}{2}=0\,.
\end{multline*}
So, we can continue our squaring process and obtain the inequality 
$$
\left(4u^2-(1-u^2)(M_-^2+M_+^2)\right)^2-4M_-^2M_+^2(1+u^2)^2\ge 0
$$
to prove, which is an inequality with a rational function $F=F(u,t,s)$ 
of $u,t,s$ on the left hand side. 

To find this rational function explicitly, one can execute the following
sequence of commands in Mathematica:

\bigskip
\begin{verbatim}
Y[u_,t_]=((1+t^2)*u+2*t)/(1-t^2);
MM[u_,t_,s_]=(Y[u,t]+Y[u,s])^2/(4+(Y[u,t]-Y[u,s])^2);
F[u_,t_,s_]=(4*u^2-(1-u^2)*(MM[u,-t,-s]+MM[u,t,s]))^2-
                          4*MM[u,-t,-s]*MM[u,t,s]*(1+u^2)^2;
Print[Factor[F[u,t,s]]];                          
\end{verbatim}

\bigskip\noindent\hglue0pt
The output looks a bit ugly with squares of two huge polynomials 
in the denominator but one can easily realize that those polynomials
come from the non-negative factors $4+(Y[u,\pm t]-Y[u,\pm s])^2$, so, executing
three more commands 

\bigskip
\begin{verbatim}
Den[u_,t_,s_]=4+(Y[u,t]-Y[u,s])^2;
G[u_,t_,s_]=Factor[F[u,t,s]*Den[u,-t,-s]^2*Den[u,t,s]^2];
Print[G[u,t,s]];                         
\end{verbatim}

\bigskip\noindent\hglue0pt
we get a much nicer output
$$
\begin{aligned}
&\frac 1{(-1+s)^8(1+s)^8(-1+t)^8(1+t)^8}
\\
&(4096(s-t)^2(s+t)^2 u^2 (1 - s t + s u - t u) (-1 + s t + s u - t u)
\\
&(3s^2-4s^4+4s^6-\dots+t^{10}u^8+s^4t^{10}u^8))
\end{aligned}
$$
(again we wrote only the very beginning and the very end of the huge polynomial
that is the most important factor).
Note that
$$
\begin{aligned}
(1-st+su-tu)(&-1+st+su-tu)=u^2(s-t)^2-(1-st)^2
\\
&\le(s-t)^2-(1-st)^2=-(1-s^2)(1-t^2)\le 0\,.
\end{aligned}
$$
So, the polynomial $P_1(u,s,t)$ to test for non-negativity
can be obtained from $G$ by executing the command

\bigskip
\begin{verbatim}
P1[u_,s_,t_]=Factor[G[u,s,t]*(1-s^2)^8*(1-t^2)^8/4096/
        (s^2-t^2)^2/u^2/(1- s*t+s*u-t*u)/(1-s*t-s*u+t*u)];                        
\end{verbatim}

\bigskip\noindent\hglue0pt
Recall that we need the non-negativity of this polynomial in the domain
$$
0\le t\le s\le \frac14\,,\qquad \frac{2s}{1+s^2}\le u\le 1\,,
$$
and our test works on $[0,1]^3$. So, we will introduce the parametrization
$$
s=z/4\,,\qquad t=yz/4\,,\qquad u=\frac{z/2}{1+z^2/16}+x
$$
and let $x,y,z$ run independently over $[0,1]$. Note that these $x,y,z$ have nothing to
do with the ones we started with. Also, parametrizing in this way, we cover
a slightly larger domain than the one we really need. With such 
parametrization, our polynomial becomes a rational function again,
so we need to multiply by the denominator, which is $(1+z^2/16)$ to some power.
To find the power, we execute the command

\bigskip
\begin{verbatim}
Print[Exponent[P1[u,s,t],u]];                        
\end{verbatim}

\bigskip\noindent\hglue0pt
which gives us $8$ as an answer.
So, the next step is to switch to our new parameters 
and to check that we, indeed, got a polynomial by executing the commands

\bigskip
\begin{verbatim}
P2[x_,y_,z_]=Factor[(1+z^2/16)^8*
                      P1[z/2/(1+z^2/16)+x,z/4,y*z/4]];
Print[P2[x,y,z]];                        
\end{verbatim}

\bigskip\noindent\hglue0pt
Looking at the output, we see that we have a huge number 
$$4722366482869645213696=4^{36}$$ in the denominator. 
That is fine because Mathematica does the computations with rational numbers 
exactly, but still we preferred to see integers only, so we executed 
one more command 

\bigskip
\begin{verbatim}
P3[x_,y_,z_]=Factor[4^36*P2[x,y,z]];                        
\end{verbatim}

\bigskip\noindent\hglue0pt
Now it is time for our test. Applied directly
to the polynomial $P_3(x,y,z)$, it reports failure, but it is enough
to replace $y$ by $1-y$ to get the ``success'' report. So, the last
two lines in our program were

\bigskip
\begin{verbatim}
R[x_,y_,z_] = P3[x,1-y,z];
LinearTest[];                       
\end{verbatim}

\bigskip
\subsection{The inequality $\dfrac{x_-}{\wt x_-}\ge \dfrac{\wt \B(x_+)-\wt \B(\wt x_+)}{\wt \B(1)-\wt \B(\wt x_+)}$}
\label{4-5}

Recall that we need to prove this inequality under the conditions
$0<x_-<x<1<x_+$. Choose $\tau$ such that $x_\pm=X(x,\pm\tau)$. We have two
restrictions on $\tau$: the condition $0<x_-<x$ implies $0<\tau<x$ and the 
condition $x_+>1$ implies $(x+\tau)^2>1-\tau^2$ or, which is the same,
$x^2+2x\tau+2\tau^2>1$.

We shall also need an explicit formula for $\wt x_-$. Solving the quadratic
equation 
$$
\frac{(\wt x_-+1)^2}{4+(\wt x_- -1)^2}=x^2\,,
$$
we get
$$
\wt x_-=\frac{2x\sqrt{2-x^2}-(1+x^2)}{1-x^2}\,.
$$
Let now 
$$
F(x,\tau)=\wt \B(X(x,\tau))=\frac 1{1+\frac{(x+\tau)^2}{1-\tau^2}}
=\frac{1-\tau^2}{1+2x\tau+x^2}\,.
$$
Since $\wt x_+=X(x,x)$ and $\wt\B(1)=\frac12$, the right hand side of 
our inequality can be rewritten as
$$
\ope{RHS}=\frac{F(x,\tau)-F(x,x)}{\frac 12-F(x,x)}\,.
$$
The left hand side is
$$
\frac{(x-\tau)(1-x^2)}{\left(2x\sqrt{2-x^2}-(1+x^2) \right)\sqrt{1-\tau^2}}\,.
$$
Since $x_-$ and the right hand side are both positive, we can rewrite our 
inequality in the form
$$
\left(2x\sqrt{2-x^2}-(1+x^2) \right)\sqrt{1-\tau^2}\le\frac{(x-\tau)(1-x^2)}{\ope{RHS}}\,.
$$
Since the expression on the right is positive, it suffices to prove the squared inequality
\begin{multline*}
\left(4x^2(2-x^2)+(1+x^2)^2-4x(1+x^2)\sqrt{2-x^2}\right)(1-\tau^2)
\\
\le\frac{(x-\tau)^2(1-x^2)^2}{\ope{RHS}^2}\,,
\end{multline*}
which is equivalent to
\begin{multline*}
4x(1+x^2)(1-\tau^2)\sqrt{2-x^2}
\\
\ge(1+10x^2-3x^4)(1-\tau^2)-\frac{(x-\tau)^2(1-x^2)^2}{\ope{RHS}^2}\,.
\end{multline*}
Since the left hand side is positive, we may square again and prove the
resulting inequality for a rational function.

All these algebraic manipulations were programmed in Mathematica as follows:

\bigskip
\begin{verbatim}
F[x_,t_]=(1-t^2)/(1+2*x*t+x^2);
RHS[x_,t_]=(F[x,t]-F[x,x])/(1/2-F[x,x]);
U[x_,t_]=4*x*(1+x^2)*(1-t^2);
V[x_,t_]=(1+10*x^2-3*x^4)*(1-t^2)-
                    (x-t)^2*(1-x^2)^2/RHS[x,t]^2;
G[x_,t_]=U[x,t]^2*(2-x^2)-V[x,t]^2;
Print[Factor[G[x,t]]];                                     
\end{verbatim}

\bigskip\noindent\hglue0pt
(here we used $t$ instead of $\tau$ and introduced two auxiliary 
functions $U$ and $V$; otherwise this program matches the above
text perfectly). The execution of this program yields the output
$$
\begin{aligned}
&-\frac 1{16(t+3x+3tx^2+x^3)^4}
\\
&\left(  (-1+x)^2(1+x)^2(-1+2t^2+2tx+x^2)(-1+5x^2)^2\right.
\\
&\left.(-1+6t^2-12t^4+\dots+150tx^{13}+25 x^{14}) \right)
\end{aligned}
$$
with some polynomial, which we will denote by $-P_1$, of quite 
reasonable size in the last parentheses.
Since $-1+2\tau^2+2\tau x+x^2>0$, we need to prove that $P_1\ge 0$.
To recover $P_1$ from $G$, we execute the command

\bigskip
\begin{verbatim}
P1[x_,t_]=Factor[G[x,t]*16*(t+3*x+3*t*x^2+x^3)^4/
              (1-x^2)^2/(-1+2*t^2+2*t*x+x^2)/(5x^2-1)^2];                                     
\end{verbatim}

\bigskip\noindent\hglue0pt
Now it is time to use our restrictions on $x$ and $\tau$.
We have \hbox{$0<\tau<x<1$} and $2\tau^2+2\tau x+x^2>1$.
The second condition is quite inconvenient to use for 
linear parametrizations, so we will replace it by a weaker
condition $\tau\ge\frac 23(1-x)$ (the left hand side is a
strictly increasing function of $\tau$ and, when 
$\tau=\frac 23(1-x)$, it equals $\frac 19(8-4x+5x^2)$, which
is less than $1$ for all $x\in(0,1)$
). Thus, we have to prove our inequality for all points 
$(x,\tau)\in\R^2$ that lie in the triangle with the vertices
$(\frac 25,\frac 25)$, $(1,0)$, and $(1,1)$.
We shall use the parametrization 
$$
(x,\tau)=(1,0)+y(0,1)+yz\left(-\tfrac 35,-\tfrac 35\right)\,,
$$
or, which is the same,
$$
x=1-\frac 35 yz\,,\qquad \tau=y-\frac 35 yz\,.
$$
When $y$ and $z$ run independently over $[0,1]$ the
point $(x,\tau)$ runs over our triangle.
This parametrization can be made by executing the command

\bigskip
\begin{verbatim}
P2[y_,z_]=Factor[5^14*P1[1-3*y*z/5, y-3*y*z/5]];                                     
\end{verbatim}

\bigskip\noindent\hglue0pt
where, again, the factor $5^{14}$ was introduced to keep all
the coefficients integer. This polynomial resisted our
attempts to prove its non-negativity by our simple test 
for several hours but finally we found the following way.
Executing the command

\bigskip
\begin{verbatim}
Print[P2[y,z]];                                     
\end{verbatim}

\bigskip\noindent\hglue0pt
and taking a quick look at $P_2$, one can see that $y$ can be factored out.
So, it is natural to divide by $y$ and introduce the polynomial $P_3$
given by 

\bigskip
\begin{verbatim}
P3[y_,z_]=Factor[P2[y,z]/y];                                     
\end{verbatim}

\bigskip\noindent\hglue0pt
An attempt to apply the linear test to $P_3$ fails too but the execution
of the commands 

\bigskip
\begin{verbatim}
R[x_,y_,z_]=Factor[2^13*P3[(1-y)/2,1-z]];
LinearTest[];
R[x_,y_,z_] = Factor[2^13*P3[1-y/2, 1-z]];
LinearTest[];                                     
\end{verbatim}

\bigskip\noindent\hglue0pt
reports success twice. Since the first pair of commands, in effect,
checks the non-negativity of $P_3$ on $[0,\frac 12]\times [0,1]$
and the second pair checks its non-negativity on 
$[\frac 12,1]\times [0,1]$, we are done.

\subsection{The inequality $M[z,y]\le M\left[\frac{z+M[z,y]}{2}, M[M[z,y],y]\right]$}
\label{4-6}

Recall that we need this inequality in the range $0\le z\le y\le 1$, $y-z\le\frac 34$.
Let $x=M[z,y]$. Let $\tau=\frac{y-z}{\sqrt{4+(y-z)^2}}$ so that
$z=X(x,-\tau)$, $y=X(x,\tau)$. Note that, since $0\le y-z\le\frac34$, we have
$0\le\tau\le\frac{y-z}2\le\frac 38<\frac{5}{13}$. Our inequality becomes 
$$
x\le M\left[\frac{X(x,-\tau)+x}{2},M[x,X(x,\tau)]\right]\,.
$$
To eliminate the square root in $X(x,\pm\tau)=\frac{x\pm\tau}{\sqrt{1-\tau^2}}$, we shall use the
substitution $\tau=\frac{2t}{1+t^2}$, $t\in[0,1)$, again. 
Note that, since the function $t\mapsto \frac{2t}{1+t^2}$ is strictly increasing
on $[0,1)$, we actually have $t\le \frac 15$.
Then 
$$
X(x,\pm\tau)=\frac{(1+t^2)x\pm 2t}{1-t^2}\,.
$$
We shall denote the the right hand side by $Y(x,\pm t)$.
Now let $U=\frac12(X(x,-\tau)+x)=\frac12(Y(x,-t)+x)$, 
$a=M[x,X(x,\tau)]=M[x,Y(x,t)]$. Note that $U$ and $a^2$ are rational
functions of $x$ and $t$. The inequality to prove is
$x\le M(U,a)$, which is equivalent to
$$
(U+a)^2\ge x^2(4+(U-a)^2)\,,
$$
or, after regrouping the terms, to
$$
2Ua(1+x^2)\ge 4x^2-(1-x^2)(U^2+a^2)\,.
$$
Since the left hand side is, clearly, non-negative, it suffices to prove
the squared inequality
$$
4U^2a^2(1+x^2)-\left(4x^2-(1-x^2)(U^2+a^2)\right)^2\ge0\,,
$$
whose left hand side is a rational function of $x$ and $t$.
To find this rational function, 
one can execute the following commands

\bigskip
\begin{verbatim}
Y[x_,t_]=((1+t^2)*x+2*t)/(1-t^2);
U[x_,t_]=(Y[x,-t]+x)/2;
Den[x_,t_]=4+(x-Y[x,t])^2;
AA[x_,t_]=(x+Y[x,t])^2/Den[x,t];
F[x_,t_]=4*U[x,t]^2*AA[x,t]*(1+x^2)^2-
              (4*x^2-(1-x^2)*(U[x,t]^2+AA[x,t]))^2;
Print[Factor[F[x,t]]];                                     
\end{verbatim}

\bigskip\noindent\hglue0pt
The output looks pretty good as is but it becomes even better if we
multiply $F$ by  $4+(x-Y(x,t))^2$, i.\,e., if we execute the commands

\bigskip
\begin{verbatim}
G[x_,t_]=F[x,t]*Den[x,t]^2;
Print[Factor[G[x,t]]];                                     
\end{verbatim}

\bigskip\noindent\hglue0pt
What we get then is 
$$
\begin{aligned}
&-\frac 1{(-1+t)^8(1+t)^8}
\\
&(16t^2(-1+tx)(-t^6+8t^3x-12t^5x+\dots-4t^6x^{10}+t^5x^{11}))
\end{aligned}
$$
with some (not really large) polynomial $P_1$ in the last parentheses.
Since $0\le tx\le 1$, we can reduce our inequality to the inequality 
$P_1(x,t)\ge 0$. To recover $P_1$ from $G$, it suffices to execute the 
command

\bigskip
\begin{verbatim}
P1[x_,t_]=Factor[G[x,t]*(1-t^2)^8/16/t^2/(1-t*x)];                                     
\end{verbatim}

\bigskip\noindent\hglue0pt
Now it is time to use the information about $x$ and $t$ we have.
Recall that $0\le t\le\frac15$. Also $x\ge\frac{2t}{1+t^2}\ge \frac32 t$
in the range of $t$ that is interesting for us and $x\le 1$. This
suggests the parametrization
$$
t=\frac z5\,,\qquad x=y+\frac 3{10}z
$$
where $y$ and $z$ run independently over $[0,1]$. The corresponding
command to execute is 

\bigskip
\begin{verbatim}
P2[y_,z_]=Factor[10^20*P1[y+3*z/10,z/5]];                                     
\end{verbatim}
 
\bigskip\noindent\hglue0pt
(we introduced the factor $10^{20}$ just to make all the coefficients
of $P_2$ integer). Now, the execution of the commands

\bigskip
\begin{verbatim}
R[x_,y_,z_]=P2[y,z];
LinearTest[];                                     
\end{verbatim}
 
\bigskip\noindent\hglue0pt
reports a successful completion of the test, thus finishing the proof.

\section{Level 5: Numerical inequalities}
 
In this section we will just prove the inequality
$$
M\left[\frac 1{2\sqrt 5},M\left[\frac 1{\sqrt 5},1\right]\right]>\frac 1{\sqrt 5}\,.
$$ 
Direct computation yields
$$
M\left[\frac 1{\sqrt 5},1\right]=\frac{\sqrt 5+1}{\sqrt{26-2\sqrt 5}}\,.
$$
We shall start with showing that the square of this number is greater than $\frac {17}{35}$.
Indeed, since $49\!>\!45$, we have $7\!>\!3\sqrt 5$, whence \hbox{$13\!>\!3(\sqrt5\!+\!2)$.} Thus, multiplying
both sides by $\sqrt 5-2>0$, we get $13(\sqrt 5-2)>3$, whence $13\sqrt 5>29$.
Now write
$$
\frac{(\sqrt 5+1)^2}{26-2\sqrt 5}=\frac{3+\sqrt 5}{13-\sqrt 5}
=\frac{39+13\sqrt 5}{169-13\sqrt 5} 
>\frac{39+29}{169-29}=\frac{68}{140}=\frac {17}{35}\,.
$$
Thus, due to monotonicity of $M$, it will suffice to prove that
$$
M\left[\frac 1{2\sqrt 5},\sqrt{\frac{17}{35}}\right]^2\ge\frac 1{ 5}\,.
$$
First, we note that $6\sqrt{119}>65$. Indeed,
$$
66-6\sqrt{119}=6(11-\sqrt{119})=\frac{12}{11+\sqrt{119}}<1\,.
$$
Now, 
\begin{multline*}
M\left[\frac 1{2\sqrt 5},\sqrt{\frac{17}{35}}\right]^2
=\frac{\left(\frac{2\sqrt{17}+\sqrt 7}{\sqrt{140}}\right)^2}
{4+\left(\frac{2\sqrt{17}-\sqrt 7}{\sqrt{140}}\right)^2}
=\frac{(2\sqrt{17}+\sqrt 7)^2}{4\cdot 140+(2\sqrt{17}-\sqrt 7)^2}
\\
=\frac{75+4\sqrt{119}}{635-4\sqrt{119}}=
\frac{225+12\sqrt{119}}{1905-12\sqrt{119}}>\frac{225+130}{1905-130}=\frac{355}{1775}=\frac 15\,,
\end{multline*}
and we are done.

The arithmetic above can be easily verified in one's head. Of course, one can ask a
computer to calculate the difference between the left and the right hand sides of
our inequality and get something like $0.000912384$, which seems to
be slightly above $0$, but this approach doesn't hold up to our declared standards of using
computers in the proofs, so, despite its shortness, we had to reject it.


\begin{thebibliography}{9999}

\bibitem{CWW}{\sc S.-Y.~A.~Chang, J.~M.~Wilson, and T.~H.~Wolff},
{\em Some weighted norm inequalities for the Schr\"odinger operator}, 
Comment. Math. Helv., {\bf60} (1985), 217--246.

\end{thebibliography}
\end{document}